\documentclass[preprint,12pt]{elsarticle}

\usepackage{amssymb}

\journal{Stochastic Processes and Applications }

\RequirePackage[OT1]{fontenc}
\RequirePackage{amsthm,amsmath, xcolor}
\RequirePackage[numbers]{natbib}
\RequirePackage[colorlinks,citecolor=blue,urlcolor=blue]{hyperref}
\usepackage{amscd,amsfonts,amssymb,latexsym,array,hhline,graphicx}
\usepackage{float}

\numberwithin{equation}{section}
\theoremstyle{plain}

\usepackage{amscd,amsfonts,amssymb,amsmath,latexsym,array,hhline,xcolor,graphicx}
\usepackage{float}

\newcommand\cE{{\cal E}}

\newcommand\cF{{\cal F}}

\newcommand\cJ{{\cal J}}

\newcommand\cL{{\cal L}}

\newcommand\cQ{{\cal Q}}

\newcommand\e{{\varepsilon}}

\def\E{{\bf E}}
\def\P{{\bf P}}
\def\Chi{{\bf 1}}

\def\bbr{{\mathbb R}}
\def\bbn{{\mathbb N}}

\def\bbz{{\mathbb Z}}

\def\text#1{\hbox{#1}}

\def\E{{\bf E}}
\def\P{{\bf P}}

\def\x{{\bf x}}
\def\Chi{{\bf 1}}

\def\d{\mathrm{d}}
\def\build #1_#2{\mathrel{\mathop{\kern 0pt #1}\limits_{#2}}} 
\newcommand{\zs}[1]{{\mathchoice{#1}{#1}{\lower.25ex\hbox{$\scriptstyle#1$}}
{\lower0.25ex\hbox{$\scriptscriptstyle#1$}}}}

\numberwithin{equation}{section}

\def\bbr{{\mathbb R}}

\def\bbr{{\mathbb R}}
\newcommand\fdem{$\Box$}
\newcommand\beq{\begin{equation}}
\newcommand\eeq{\end{equation}}
\newcommand\bea{\begin{eqnarray}}
\newcommand\eea{\end{eqnarray}}
\newcommand\bean{\begin{eqnarray*}}
\newcommand\eean{\end{eqnarray*}}

\newtheorem{theo}{Theorem}[section]
\newtheorem{prop}[theo]{Proposition}
\newtheorem{lemm}[theo]{Lemma}

\newtheorem{coro}[theo]{Corollary}
\newtheorem{rem}[theo]{Remark}

\begin{document}

\begin{frontmatter}



\title{The ruin problem for L\'evy-driven linear stochastic equations with applications to actuarial models with negative risk sums}

\author[A1]{Yuri KABANOV } 
\author[A2]{Serguei  PERGAMENSHCHIKOV}
\address[A1]{Universit\'e Bourgogne Franche-Comt\'e, Laboratoire de Math\'ematiques,\\  16 Route de Gray,  25030 Besan\c{c}on cedex, France, and\\ 
National Research University - Higher School of Economics,  Moscow, Russia\\
 Email: youri.kabanov@univ-fcomte.fr}

\address[A2]{Universit\'e de Rouen Normandie, Laboratoire de Math\'ematiques Rapha\"el Salem, 
Technop\^ole du Madrillet,
76801 Saint-\'Etienne-du-Rouvray, France, and\\
National Research Tomsk State University, 
International Laboratory of Statistics of Stochastic Processes and Quantitative Finance, Tomsk, Russia\\
 Email: Serge.Pergamenchtchikov@univ-rouen.fr}

\begin{abstract}
 We study the asymptotic of the ruin probability for a process which is the  solution of  linear SDE defined by a pair of independent L\'evy processes.
Our main interest is the  model describing the evolution of the capital reserve of an insurance company selling annuities and investing in a risky asset.  Let $\beta>0$ 
be the root of the cumulant-generating function $H$ of the  increment of the log price  process $V$. 
We show that  the ruin probability admits the exact asymptotic $Cu^{-\beta}$ as the initial capital $u\to\infty$  assuming only that   the law  of $V_T$ is non-arithmetic without   any further assumptions on the price process. 
\end{abstract}



\end{frontmatter}

 \section{Introduction} 
 
The general ruin problem can be formulated as follows.  We are given a family of scalar processes  $X^u$ with the initial values $u> 0$.  The object of interest is  the  
 exit probability of $X^u$ from the positive half-line as a function of $u$. More formally, let $\tau^u:=\inf \{t\colon X^u_t\le 0\}$.   
 The question is to determine the function $\Psi (u,T):=\P(\tau^u\le T)$ (the ruin probability  on a finite interval $[0,T]$) or  $\Psi (u):=\P(\tau^u< \infty)$ (the ruin probability  on $[0,\infty[$).  The exact solution of the problem is available only in  rare cases.  For instance, for  $X^u=u+W$ where $W$ is the Wiener process we have  $\Psi (u,T)=\P(\sup_{t\le T}W_t\ge u)$  and it remains to recall that the explicit formula for the distribution of the supremum of the Wiener process was obtained already in the Louis Bachelier thesis of 1900 which is, probably, the first ever mathematical study on continuous stochastic processes. Another example is the well-known explicit formula for $\Psi (u)$ in the Lundberg model of the ruin of insurance company with exponential claims.  Of course,  for more complicated cases the explicit formulae are not available  and only asymptotic results or bounds can be obtained as it is done, e.g., in the Lundberg--Cram\'er theory. 

 In this paper we consider the ruin problem  for a rather general model, suggested by Paulsen in  \cite{Paul-93},  in which 
$X^u$ is given as the  solution of   linear stochastic  equation (sometimes called the generalized Ornstein--Uhlenbeck process)
 \beq
 \label{basic}
 X^u_t=u+P_t+\int_{]0,t]}X^u_{s-}dR_s,
 \eeq
 where  $P$ and $R$ are independent L\'evy processes   with the 
 L\'evy triplets  $(a,\sigma^2,\Pi)$ and $(a_P,\sigma_P^2,\Pi_P)$, respectively.  {\bf We assume} that $\Pi (]-\infty,-1])=0$ (otherwise $\Psi(u)=1$ for all $u>0$) and  $P$ is a not a subordinator (otherwise     $\Psi(u)=0$ for all $u>0$ since the process  $X^u$ is strictly positive, see (\ref{uY}), (\ref{Y_t})).  Also {\bf we exclude}  the case  $R=0$ well studied in the literature, see \cite{KKM}.

 There is a growing  interest in  models of this type because they describe the evolution of reserves of insurance companies investing in a risky asset with the price process $S$.  In the actuarial context $R$ is interpreted as the relative price  
  process with $dR_t=dS_t/S_{t-}$, that is 
the price process $S$ is  the stochastic (Dol\'eans) exponential $\cE(R)$. 
The log price process  $V=\ln \cE(R)$ is also a L\'evy process with the triplet $(a_V,\sigma^2,\Pi_V)$. 
Recall that the behavior of the ruin probability in such models  is radically different from that in the classical actuarial  models. 
For instance, if the price of the risky asset follows a geometric Brownian motion, that is, $R_t=at+\sigma W_t$, and the risk process $P$ is as in the Lundberg model, then  
$\Psi(u)=O(u^{1-2a/\sigma^2})$, $u\to \infty$, 
 if $2a/\sigma^2> 1$,  and 
$\Psi(u)\equiv 1$ otherwise,   \cite{FrKP, PeZe-06, KP}.

We are especially interested in the case where the process $P$ describing the  ``business 
part" of the model has only upward jumps (in other words, $P$ is spectrally positive). In  the classical actuarial literature such models are  referred to as the annuity insurance models (or models with negative risk sums), \cite{Gr, Sax}, while in modern sources  they serve also to describe   the capital reserve  of a venture company investing in development of new technologies  and selling innovations; sometimes they are referred to as the dual models,  \cite{ABL, Asm, AGS, Bayraktar-Egami}, etc.  

The mentioned specificity of models with negative risk sums leads to a continuous  downcrossing of the zero level by the capital reserve process. This allows us to obtain the exact (up to a multiplicative constant)  asymptotic of the ruin probability under weak assumptions on the price dynamics.

Let $H: q\mapsto \ln \E\,e^{-q V_1}$ be the  cumulant-generating  function  of the increment of  log price process $V$ on the interval $[0,1]$. 
The function $H$ is convex and its effective domain ${\rm dom}\,H$ is a convex subset of $\bbr$ containing zero. 

It is well-known that the asymptotic of the ruin probability $\Psi(u)$ as $u\to \infty$  is determined by the strictly positive root $\beta$ of $H$, assumed existing  and  laying in the interior of ${\rm dom}\,H$. Unfortunately, the existing results are overloaded by numerous integrability assumptions on processes $R$ and $P$ while the law $\cL(V_T)$ of the random variable $V_T$ is 
required to contain an absolute continuous component where $T$ is independent random variable uniformly distributed on $[0,1]$, see, e.g., Th. 3.2 in \cite{Paul-02}.

The aim of our study  is to obtain the  exact asymptotic of the exit probability under  the weakest conditions. 
Our main result has the following easy to memorize formulation. 

 \begin{theo}
\label{Main2+}
Suppose that $H$ has a root  $\beta>0$  laying in  ${\rm int}\,{\rm dom}\,H$ and 
$\int |x|^{\beta}I_{\{|x|>1\}}\Pi_P(dx)<\infty$. 
Then 
\beq
0< \liminf_{u\to \infty} u^{\beta}\Psi(u)\le \limsup_{u\to \infty} u^\beta \Psi(u)<\infty.  
\eeq

If, moreover, $P$ has only upward jumps and the distribution  
$\cL(V_1)$  is non-arithmetic, then $\Psi(u)\sim C_{\infty}u^{-\beta}$ where $C_{\infty}>0$. 
\end{theo}




In  our argument  we are based on the theory of distributional equations as presented in the paper by  Goldie, \cite{Go91} 
and  on the criterion by Guivarc'h and Le Page, \cite{GL}, which simple proof can be found in the recent paper \cite{BD} by Buraczewski  and Damek. This criterion gives a necessary and sufficient condition for the strict positivity of the constant in the Kesten--Goldie theorem determining the rate of decay of the tail of solution at infinity. Its obvious corollary   allows us to simplify radically the proofs  and get rid of additional  assumptions presented in the earlier papers, see \cite{Paul-93, Paul-98, Paul-02, Nyrh-99, Nyrh-01, Kalash-Nor, BKM}  and references therein.  Our technique involves only affine distributional equations and avoids more demanding Letac-type equations.  

\smallskip

%

The structure of the paper is the following. In Section \ref{prel} we formulate the model and provide some prerequisites from L\'evy processes.  
Section \ref{reduc}  contains a well-know reduction of the ruin problem to the  study of 
asymptotic behavior of a stochastic integral.  
In Section \ref{moments} we prove moment inequalities for maximal functions of stochastic integrals needed to analysis of the limiting  behavior of an exponential functional in Section \ref{asymp}. The latter section is concluded by the proof of the main result and some comments on its formulation.  
In Sec\-tion \ref{prob1} we establish Theorem \ref{ruintheo} on the ruin with probability one using the technique suggested in \cite{PeZe-06}. This theorem implies, in  particular, that  in the classical  model with negative risk sums and investments 
in the risky asset with  price following a geometric Brownian motion the ruin is imminent if $a\le \sigma^2/2$, \cite{KP}. 
In Section \ref{Sec6} we discuss examples.  Our presentation is oriented towards  the reader  with preferences in the L\'evy processes rather than in the theory 
of distributional equations (called also implicit renewal theory). That is why in Section \ref{sec:A} (Appendix) we provide a rather detailed information on the latter
covering the arythmetic case. In particular, we give a proof of a version   
of Grincevi\u{c}ius theorem under slightly weaker conditions as in the original paper. 

\smallskip 
 We express our gratitude to the anonymous referees whose constructive criticism 
lead us to substantial strengthening of the main result as well as to E. Damek, D. Buraczewski, and Z. Palmowski   who communicated to the authors a number of useful references on distributional equations.

\section{Preliminaries from  the theory of L\'evy processes}
\label{prel}

Let  $(a,\sigma^2,\Pi)$ and $(a_P,\sigma_P^2,\Pi_P)$ be the L\'evy triplets 
of the processes $R$ and $P$ corresponding to the standard\footnote{Other truncation functions are also used in the literature, see, e.g., \cite{Paul-02}} truncation function
$h(x):=xI_{\{|x|\le 1\}}$.  Putting $\bar h(x):=xI_{\{|x|> 1\}}$ we can write 
the canonical decomposition of $R$ in the form
\beq
\label{Rt}
R_t=at+\sigma W_t+h*(\mu -\nu)_t+\bar h*\mu_t
\eeq
where $W$ is a standard Wiener process,  the Poisson random measure $\mu(dt,dx)$ is the jump measure of $R$ having the (deterministic) compensator 
$\nu(dt,dx)=dt\Pi (dx)$.  For  notions and results see \cite{JS}, Ch. 2 and also \cite{CT}, Chs. 2 and 3.

\smallskip
As in \cite{JS}, we use $*$ for the standard notation of stochastic calculus  for integrals with respect to random measures. For instance,    
$$
h*(\mu -\nu)_t=\int_0^t\int h(x)(\mu-\nu)(ds,dx). 
$$
We hope that the reader will be not confused that $f(x)$  may denote the whole function $f$ or its value at $x$; the typical example is $\ln (1+x)$ explaining why such a flexibility is convenient. The symbols $\Pi(f)$ or  $\Pi(f(x))$ stands for the integral of $f$ with respect to the measure $\Pi$. 

Recall that
 $$
\Pi(|x|^2\wedge 1):=\int (|x|^2\wedge 1) \Pi(dx)<\infty
 $$
and the condition $\sigma=0$ and $\Pi(|h|)<\infty$ is necessary and sufficient for  $R$ to have trajectories of (locally) finite variation, see Prop. 3.9 in \cite{CT}.  

The process $P$ describing the actuarial (``business") part of the model admits a similar representation:  
\beq
\label{Pt}
P_t=a_Pt+\sigma_P W^P_t+h*(\mu^P -\nu^P)_t+\bar h*\mu^P_t. 
\eeq

The L\'evy processes $R$ and $P$ generate  the filtration ${\bf F}^{R,P}=(\cF^{R,P}_t)_{t\ge 0}$.  
\medskip

\smallskip 
\noindent
{\bf Standing assumption S.0} 
 {\em The L\'evy measure $\Pi$ is concentrated on the interval $]-1,\infty[$; $\sigma^2$ and $\Pi$ do not vanish simultaneously; the process $P$ is not a subordinator. }

\medskip

Recall that  subordinator is an increasing L\'evy process. Accordingly to \cite{CT}, Prop. 3.10,  
the process $P$ is not a subordinator  if and only if $\sigma^2_P>0$, or 
one of the following three conditions hold: 

1) $\Pi_P(]-\infty,0[)>0$, 

2) $\Pi_P(]-\infty,0[)=0$, $\Pi_P(xI_{\{x>0\}})=\infty$, 

3) $\Pi_P(]-\infty,0[)=0$, $\Pi_P(xI_{\{x>0\}})<\infty$, 
$\Pi_P(xI_{\{x>0\}})-a_P>0$. 
 
\medskip
 
In the context of financial models the stochastic exponential 
$$
\cE_t(R)=e^{R_t -\frac 12\sigma^2t+\sum_{s\le t}(\ln (1+\Delta R_s)-\Delta R_s)}
$$
stands for the price of a risky asset (e.g., stock).  The log price $V:=\ln \cE(R)$ is a 
L\'evy process and can be written in the form 
\beq
\label{V}
V_t=at-\frac 12 \sigma^2 t + \sigma W_t+ h*(\mu-\nu)_t+(\ln (1+x)-h)*\mu_t. 
\eeq
Its L\'evy triplet is  $(a_V,\sigma^2,\Pi_V)$ where 
$$
a_V=a-\frac{\sigma^2}{2}+\Pi(h(\ln(1+x))-h)  
$$ and $\Pi_V=\Pi\varphi^{-1}$,  $\varphi: x\mapsto \ln (1+x)$.  

\smallskip
The cumulate-generating function $H:q\to \ln  \E\,e^{-qV_1}$ of the random variable $V_1$ admits an explicit expression. Namely, 
\beq
\label{sec.Main.1}
H(q):=-a_V q+\frac{\sigma^{2}}{2}q^{2}
+\Pi \big(e^{-q \ln(1+x)}-1+q h(\ln (1+x))\big).  
\eeq
Its effective domain ${\rm dom}\,H=\{q\colon\ H(q)<\infty\}$ is the set $\{J(q)<\infty\}$ where 
\begin{equation}
\label{sec.Main.1-00}
 J(q):=
\Pi \big(I_{\{|\ln (1+x)|>1\}}\,e^{-q\ln (1+x)}\big)=\Pi \big(I_{\{|\ln (1+x)|>1\}}(1+x)^{-q}\big). 
\end{equation}
Its interior is the open interval  $]\underline q,\bar q[$ with 
$$
\underline q:=  \inf\{q\le 0\colon\, J(q)<\infty\},\qquad   \bar q:=\sup\{q\ge 0\colon\, J(q)<\infty\}.
$$
Being a   convex function, $H$  is continuous and 
admits finite right and left derivatives on  $]\underline q,\bar q[$.  If $\bar q>0$, then the right derivative 
$$
D^+H(0)=-a_V-\Pi(\bar h(\ln(1+x)))<\infty,
$$ 
though it may be equal to $-\infty $. 

In formulations of our asymptotic results  we shall always assume that $\bar q>0$ and the equation $H(q)=0$ has  a root $\beta\in ]0,\bar q[$. Since $H$
is not constant, such a root is unique.  Clearly, it exists if and only if  $D^+H(0)<0$ and 
$\limsup_{q\uparrow\bar q}H(q)/q>0$. In the case where  $\underline q<0$  the condition  $D^-H(0)>0$
is necessary to ensure that $H(q)<0$ for sufficiently small in absolute value $q<0$.

\smallskip
 
 If $J(q)<\infty$, then   
 the process $m=(m_t(q))_{t\le 1}$ with 
\beq
\label{mart}
m_t(q):=e^{-qV_{t}-tH(q)}
\eeq
is a martingale and
\begin{equation}
\label{sec.Main.1-0}
 \E\,e^{-qV_{t}}= e^{tH(q)}, \qquad t\in [0,1].      
\end{equation}
In particular, we have that $H(q)=\ln \E\,e^{-qV_{1}}=\ln \E M^q$. For the above properties see, e.g., Th. 25.17 in \cite{Sato1999}.

\smallskip
Note that    
\beq
\label{martbound}
\E\,\sup_{t\le 1}e^{-qV_{t}}<\infty \qquad \forall\,q\in ]\underline q,\bar q[.    
\eeq  
 Indeed, let $q\in ]0,\bar q[$. 
Take $r\in ]1,\bar q/q[$. Then $\E\,m^r_1(q)=e^{H(qr)-rH(q)}<\infty$. By virtue of the Doob inequality the maximal function 
$m^*_1(q):=\sup_{t\le 1}m_t(q)$ belongs to $L^r$  and it  remains to observe  that 
$e^{-qV_{t}}\le C_q m_t(q)$ where the constant $C_q=\sup_{t\le 1}e^{t H(q)}$. Similar arguments work for $q\in]\underline q,0[$.

\section{Ruin problem: a reduction}
\label{reduc}

Let us introduce the process 
\beq
\label{Y_t}
Y_t:=-\int_{]0,t]} \cE^{-1}_{s-}(R)dP_s=-\int_{]0,t]}  e^{-V_{s-}}dP_s.  
\eeq
Due to independence of $P$ and $R$ the joint quadratic characteristic $[P,R]$ is zero, 
and the straitforward application  of the product formula for semimaringales 
shows that the process 
\beq
\label{uY}
X_t^u:=\cE_t(R)(u-Y_t)
\eeq
 solves the non-homogeneous linear  equation (\ref{basic}), i.e. 
 the solution of the latter is given by this stochastic version of the Cauchy formula.

The positivity of $\cE(R)$ implies that    $ \tau^u=\inf \{t\ge 0:\ Y_{t} \ge u\}$.


\smallskip 

The following lemma is due to Paulsen and Gjessing, see \cite{Paulsen1993}.  
\begin{lemm}
\label{G-Paulsen} 
If $Y_t\to Y_{\infty}$ almost surely as $t\to \infty$  where  $Y_{\infty}$ is a finite random variable unbounded from above, 
then  for all  $u>0$  
\beq
\label{Paulsen}
\bar G(u)\le\, \Psi(u)=\frac{\bar G(u)}{\E\left(\bar G(X_{\tau^{u}})\, \vert\, \tau^{u}<\infty\right)}\le \frac{\bar G(u)}{\bar G(0)},
\eeq  
where $\bar G(u):=\P(Y_\infty>u)$. 

In particular, if $\Pi_P(]-\infty,0])=0$, then $\Psi(u)=
\bar G(u)/\bar G(0)$.
\end{lemm}
\noindent
{\sl Proof.}
Let $\tau$ be an arbitrary stopping time with respect to the  filtration $(\cF^{P,R}_t)$. 
As we assume that the finite limit $Y_\infty$ exists,  the random variable 
$$
Y_{\tau,\infty}:=\begin{cases}
-\lim_{N\to \infty } 
\int_{]\tau,\tau+N]}\,e^{-(V_{t-}-V_{\tau})}dP_{t},&\tau<\infty, \\
0, & \tau=\infty,
\end{cases} 
$$ 
is well defined.  On the set $\{\tau<\infty\}$
\beq
\label{YX}
Y_{\tau,\infty}=e^{V_\tau}(Y_{\infty}-Y_\tau)=X_{\tau} +e^{V_\tau}(Y_\infty-u). 
\eeq   Let $\xi$ be a $\cF_{\tau}^{P,R}$-measurable random variable.  Since the L\'evy process $Y$ starts afresh at $\tau$,  the conditional distribution of $Y_{\tau,\infty}$ given $(\tau,\xi)=(t,x)\in \bbr_+\times\bbr$ is the same as the distribution of $Y_{\infty}$. It follows that   
$$
\P\left(
Y_{\tau,\infty}>\xi, \ 
\tau<\infty
\right)
=\E\,\bar G(\xi)\, \Chi_{\{ \tau<\infty\}}.
$$
Thus, if $\P(\tau<\infty)>0$, then
$$
\P\left(Y_{\tau,\infty}>\xi, \ 
\tau<\infty
\right)
=\E\left(\bar G(\xi)\, \vert\, \tau<\infty\right)\,\P(\tau<\infty)\,.
$$
Noting that  $\Psi(u):=\P(\tau^{u}<\infty)\ge \P(Y_{\infty}>u)>0$,   we deduce from here using (\ref{YX}) that 
\begin{align*}
\bar G(u)&=\P\left(
Y_{\infty}>u,\ \tau^{u}<\infty\right)=
\P\left(Y_{\tau^u,\infty}>X_{\tau^{u}},\ 
\tau^{u}<\infty
\right)\\
&=\E\left(\bar G(X_{\tau^{u}})\, \vert\, \tau^{u}<\infty\right)\,\P(\tau^{u}<\infty)
\end{align*}
implying the equality in (\ref{Paulsen}). The result follows since  
$X_{\tau^{u}}\le 0$ on the set  $\{\tau^u<\infty\}$ and,  in the case  
 where  $\Pi_P(]-\infty,0])=0$, the process $X^u$ crosses zero in a continuous way, i.e.  $X_{\tau^{u}}= 0$ on this set. \fdem 
 
In view of the above lemma the proof of  Theorem \ref{Main2+} is reduced to establishing of   the existence of finite limit $Y_{\infty}$ and finding the asymptotic of the tail of its distribution.

\section{Moments of the maximal function} 
\label{moments}
In this section we prove a simple but important result on the existence of moments of the maximal function of the process $Y$ on the interval $[o,1]$, i.e. of the random variable $Y_1^{*}:=\sup_{t\le 1} |Y_t|$. 

Before the formulation we recall  the Novikov  inequalities, \cite{Novikov}, also referred to as the Bichteler--Jacod inequalities, see  \cite{BJ, MR},  providing bounds  for   moments of the maximal function $I^*_1$ of stochastic integral $I=g*(\mu^P-\nu^P)$ where $g^2*\nu^P_1<\infty$.    In dependence of the  parameter $\alpha\in [1,2]$ they have the following form: 
\beq
\label{Novikov}
\E I_1^{*p}\le C_{p,\alpha}
\begin{cases} \E\,\big (|g|^\alpha*\nu^P_1\big)^{p/\alpha}, \quad  &\forall\, p\in ]0,\alpha],\\ 
\E\,\big (|g|^\alpha*\nu^P_1\big)^{p/\alpha}+\E\, |g|^p*\nu^P_1, \quad & \forall\, p\in [\alpha,\infty[.
\end{cases}
\eeq
\begin{lemm}\label{Q}
Let $p>0$ be such that $\Pi_P(|\bar h|^{p})+\E\,\sup_{t\le 1}e^{-pV_{t}}<\infty$. Then  
$\E\,Y_1^{*p}<\infty$.  
\end{lemm}

\noindent 
{\sl Proof.}  We start with the case where $p\in ]0,1[$.  The elementary inequality $(\sum x_k)^p\le \sum x_k^p$  allows  us to treat separately the  integrals corresponding to each  term in the representation 
$$
P_t=a_Pt+ \sigma_PW^P_t+h*(\mu^P -\nu^P)_t+\bar h*\mu^P_t. 
$$
 Recall that in the detailed notations 
$f*\mu^P_1=\sum_{\{s\le 1:\  \Delta P_s>0\}}f(s,\Delta P_s)$ and $V_-=(V_{s-})$.   Using the mentioned  inequality  we get that  
\bea
\nonumber
\E\,(e^{-V_{-}}|\bar h|*\mu^P_1)^{p}&\le& \E\, e^{-{p V_{-}}}|\bar h|^{p}*\mu^P_1=\Pi_P(|\bar h|^{p}) \E\int_0^1
e^{-{p V_{t}}}dt\\
\label{1}
&\le& \Pi_P(|\bar h|^{p}) \E\,\sup_{t\le 1}e^{-pV_{t}}. 
\eea
Note that 
\beq
\label{2}
\E\,\Big(\int _0^1e^{-V_{t}}dt \Big)^{p}\le \E\,\sup_{t\le 1}e^{-pV_{t}}.
\eeq
By the Burkholder--Davis--Gundy  inequality 
\beq
\label{3}
\E\sup_{t\le 1}\Big\vert \int _0^te^{-V_{s}}dW^P_s \Big\vert^{p}\le
C_{p}\E\Big (\int _0^1e^{-2V_{s}}ds\Big)^{p/2}\le 
C_{p}\E \sup_{t\le 1} e^{-pV_{t}}. 
\eeq
Using the Novikov inequality  (with $\alpha=2$) we have 
\bea
\nonumber
\E\sup_{t\le 1}\big\vert e^{-V_-}h*(\mu^P -\nu^P)_t \big\vert^{p}
&\le&
C_{p,2}(\Pi(h^2))^{p/2}\E\Big (\int_0^1e^{-2V_{t}}dt\Big)^{p/2}\\
\label{4}
&\le& 
C_{p,2}(\Pi(h^2))^{p/2} \E\, \sup_{t\le 1}e^{-pV_{t}}. 
\eea
\smallskip
From  these estimates and the property (\ref{martbound}) we have that  $\E\, Y_*^p<\infty$.

Let $p\in ]1,2[$. 
By the Novikov inequality with $\alpha=1$ and  we have: 
\bean
\E\,\sup_{t\le 1}|e^{-V_-}\bar h*(\mu^P-\nu^P)_t|^{p}
&\le&
 C_{p,1}\Big( \big(\E\, (e^{-V_-}\bar h*\nu_1^P\big)^{p} +\E\, e^{-{p}V_-}\bar h^{p}*\nu_1^P\Big) \\
&\le&  \tilde C_{p,1} \E\, \sup_{t\le 1}e^{-pV_t}<\infty,
\eean
where $\tilde C_{p,1}:=C_{p,1}\big(\big(\Pi_P(\bar h) \big)^{p}+\Pi_P(\bar h^{p})\big)$. 
Using again the Novikov inequality but with $\alpha=2$ we obtain that   
\begin{align*}
\E\,\sup_{t\le 1}|e^{-V_-}h*(\mu^P-\nu^P)_t|^{p}&\le C_{p,2}\E\, (e^{-2V_-}h^2*\nu^P_{1})^{p/2}\\
&\le  C_{p,2}(\Pi_P(h^2))^{p} \E\, \sup_{t\le 1}e^{-pV_t}
<\infty.
\end{align*}
Estimates for the integrals with respect to  $dt$ and $dW^P$ are of the same form as for the previous case. Using  the inequality for the $L^p$-norm of the sum, we get that  $\E\,Y^{*p}_1<\infty$.  
 
 Finally, let $p\ge 2$. 
Using  the Novikov  inequality  with $\alpha = 2$,  we have:
\bean
\E \sup_{t\le 1}|g*(\mu^P-\nu^P)_t|^{p}&\le& C_{p,2}
\big(\Pi_P (|x|^{2})\big)^{p/2}\, \E\Big(\int^{1}_{0}\,e^{-2V_{t}}d t \Big)^{p/2}\\
&&
+C_{p,2} \Pi_P (|x|^{p})\,\E\int^{1}_{0}\,e^{-pV_{t}}d t\\
&\le & C_{p,2}  \big(
\big(\Pi_P (|x|^{2})\big)^{p/2}
+
\Pi_P (|x|^{p})
\big)
\E\, \sup_{t\le 1}e^{-pV_t}<\infty.
\eean
Again the arguments for the integrals with respect to $dt$ and $dW^P$ remain valid.  \fdem

\section{Convergence of $Y_t$} 
\label{asymp}

Using Lemma \ref{Q} the convergence  $Y_t$ as $t\to \infty$  can be easily established under very weak assumptions. Namely, we have the following:  

\begin{prop} 
\label{YYY}
If there is $p>0$ such that $H(p)<0$, and  $\Pi_P(|\bar h|^p)<\infty$,    then $Y_t$ converge a.s. to a finite random variable $Y_{\infty}$ unbounded from above and solving the distributional equation 
\begin{equation}
\label{deq}
Y_{\infty} \stackrel{d}{=}Y_1+M_1\,Y_{\infty}\,,\quad
Y_{\infty} \ \ \mbox{independent of}\ (M_1,Y_1),
\end{equation}
where $M_1:=e^{-V_1}$. 

\end{prop}
\noindent
{\sl Proof.} Without loss of generality we  assume that $p<1$ and $H(p+)\neq \infty$.   
For any integer  $j\ge 1$ 
 $$
Y_j-Y_{j-1} =
 M_1 \dots M_{j-1}Q_j,  
\qquad.  
$$
where   $(M_j,Q_j)$ are independent random variables,  
\beq
\label{MQ}
M_j:=e^{-(V_{j}-V_{{j-1}})}, \qquad Q_j:=-\int_{]{j-1},j]}\,e^{-(V_{v-}-V_{{j-1}})}dP_v 
\eeq
with distributions  $\cL(M_j)=\cL(M_1)$ and $\cL(Q_j)=\cL(Y_1)$.  By  assumption, $\rho:=\E M_1^p=e^{H(p)}<1$  and $\E |Y_1|^p<\infty$ in virtue  of (\ref{martbound}) and Lemma \ref{Q}. 
Since   
$\E M_1...M_{j-1}|Q_j|=\rho^j\E |Y_1|^p$ we have that 
$\E \sum_{j\ge 1 } |Y_j-Y_{j-1}|^p<\infty$ and, hence, $\sum_{j\ge 1 } |Y_j-Y_{j-1}|^p<\infty$ a.s. 
 But then also $\sum_{j\ge 1} |Y_j-Y_{j-1}|<\infty$ a.s.
and, therefore, the sequence $Y_n$ converges almost surely to some finite random variable $Y_\infty$.

Put
$$
\Delta_{n}:=\sup_{n-1\le v\le n}
\left\vert 
\int_{]n-1,v]}\,e^{-V_{s-}}\,d P_{s}
\right\vert , \qquad n\ge 1.
$$
Note that  
$$
\E\,\Delta_{n}^{p}
=
\E\,\prod^{n-1}_{j=1}\,M^p_{j}\,
\sup_{n-1\le v\le n}
\left\vert 
\int_{]n-1,v]}\,e^{-(V_{s-}-V_{n-1})}\,d P_{s}
\right\vert^{p} 
=\rho^{n-1}\,\E\,Y_1^{*p}\,<\infty. 
$$
For any $\e>0$ 
we get using the Chebyshev inequality that 
$$
\sum_{n\ge 1}\,\P(\Delta_{n}>\e)\le \e^{-p}\E\,Y_1^{*p}\sum_{n\ge 1}\rho^{n-1}<\infty. 
$$

By the Borel--Cantelli  lemma  $\Delta_n(\omega)\le \e$ for all  $n\ge n_0(\omega)$ for each $\omega\in \Omega$ except a null-set.  This implies the convergence  $Y_t\to Y_{\infty}$ a.s., $t\to \infty$.

\smallskip
Let us consider the sequence 
$$
Y_{1,n}:=Q_2+M_2Q_3+\dots+M_2\dots M_{n}Q_{n+1}
$$
converging almost surely to a  random variable   $Y_{1,\infty}$ distributed as $Y_{\infty}$. 
Passing to the limit in the obvious  identity
$Y_n=Q_1+M_1Y_{1,n-1}$
 we  obtain that $Y_\infty=Q_1+M_1Y_{1,\infty}$. 
 For finite $n$ the random variables $Y_{1,n}$ 
and $(M_1,Q_1)$ are independent, $\cL(Y_{1,n})=\cL(Y_n)$. 
Hence, $Y_{1,\infty}$ 
and $(M_1,Q_1)$ are independent, $\cL(Y_{1,\infty})=\cL(Y_\infty)$ and  
 $\cL(Y_{\infty})=\cL(Q_1+M_1Y_{1,\infty})$.  This is exactly the properties abbreviated by (\ref{deq}). 

\smallskip
It remains to check that $Y_{\infty}$ is unbounded from above. For this  
it is useful the following simple observation. 

\begin{lemm} If the random variables $Q_1$ and $Q_1/M_1$ are unbounded from above, then  
 $Y_{\infty}$ is also unbounded from above. 
\end{lemm}
\noindent
{\sl Proof.}   Since 
$Q_1/M_1$ is unbounded from above,  we have,  due to independence of $(Q_1/M_1)$  and $Y_{1,\infty}$,  that $\P(Y_{1,\infty}>0)=\P(Y_{\infty}>0)>0$. Take arbitrary $u>0$.  Then 
\bean
\P(Y_{\infty}> u)&\ge& \P(Q_1+M_1Y_{1,\infty}>u,\, Y_{1,\infty}> 0)\ge \P(Q_1>u,\, Y_{1,\infty}> 0)\\
&=& \P(Q_1>u)\P(Y_{1,\infty}> 0)>0
\eean
and the lemma is proven. \fdem

\medskip

\noindent
{\bf Notation:}   ${\cJ}_\theta:=\int^{1}_{0}\, e^{-\theta V_{v}}d v$, \  ${Q}_\theta:=-\int^{1}_{0}\, e^{-\theta V_{v-}}dP_v$ where $\theta=1$ or $-1$.  
\smallskip

\begin{lemm} 
\label{Q/M}
$\cL(Q_{-1})=\cL(Q_1/M_1)$.
\end{lemm}
\noindent
{\sl Proof.} We have: 
\bean
\int_{]0,1]} \sum_{k=1}^n e^{V_{k/n-}}I_{](k-1)/n,k/n]}(v)dP_v&= &\sum_{k=1}^n e^{V_{k/n}}
(P_{k/n}-P_{(k-1)/n}),\\
e^{V_1}\int_{]0,1]} \sum_{k=1}^n e^{-V_{k/n-}}I_{](k-1)/n,k/n]}(v)dP_v&= &\sum_{k=1}^n e^{V_1-V_{k/n}}
(P_{k/n}-P_{(k-1)/n}).
\eean
Note that $V$ and $P$ are independent, the increments  $P_{k/n}-P_{(k-1)/n}$ are independent and identically distributed, and $\cL(V_1-V_{k/n})=\cL(V_{(n-k)/n})$. 
Thus, the right-hand sides of the above identities have the same distribution. The result follows  because the left-hand sides  tend in probability, respectively, to  $-Q_{-1}$ and $-Q_1/M_1$.  \fdem

\medskip Thus, $Y_\infty$ is unbounded from above if so are the stochastic integrals $Q_\theta$. Lemma  \ref{suppQ} below shows that $Q_\theta$ are unbounded from above if the  ordinary integrals ${\cJ}_\theta$ are unbounded from above. For the latter property we prove  necessary and sufficient conditions in terms of defining characteristics (Lemma \ref{suppQ.1}). The case where these conditions are not fulfilled we treat separately (Lemma \ref{twocases}).

\medskip

\begin{lemm}
\label{suppQ} 
If   $\cJ_\theta$ is unbounded from above, so is  $Q_\theta$. 
\end{lemm} 

\noindent 
{\sl Proof.}
We argue using the following  observation: if $f(x,y)$ is a measurable function and $\xi$, $\eta$ are independent random variables with distributions $\P_\xi$ and $\P_\eta$, then  the distribution of $f(\xi,\eta)$ is unbounded from below if the distribution of  $f(\xi,y)$ is unbounded from below on a set of $y$ of positive measure $\P_\eta$.  

In the case $\sigma_P^2>0$, we use the representation  
$$
Q_\theta=-\sigma_P\int_{]0,1]} e^{-\theta V_{v-}}dW^P_v + \int_{]0,1]} e^{-\theta V_{v-}}d(\sigma _PW^P_v-P_v). 
$$   
Applying the above observation (with $\xi=W^P$ and $\eta=(R,P-\sigma_P W^P)$) and taking into account that the Wiener integral of a strictly positive deterministic function is a nonzero Gaussian random variable, we get that $Q_\theta$ is unbounded. 

 Consider the case where $\sigma_P^2=0$. 

For $\e>0$ we denote by $\zeta^\e$
the locally square integrable martingale  with 
\beq
\label{zeta}
\zeta^\e_t:=\,e^{-\theta V_{-}}\,I_{\{|x|\le\e\}} x*(\mu^P-\nu^P)_{t}.  
\eeq 
Since $\langle \zeta^\e\rangle_1=e^{-2\theta V_{-}}\,I_{\{|x|\le\e\}} x^2*\nu^P_{1}\to 0$ as  $\e\to 0$, we have that  
$\sup_{t\le 1}|\zeta^\e_t|\to 0$ in probability. 

Note that 
$$
Q_\theta=(\Pi_P(xI_{\{\e\le |x|\le 1\}})-a_P) \cJ_\theta - \zeta^\e_t - e^{-\theta V_{-}}\,I_{\{|x|< \e\}} x*\mu^P_1.
$$

Take    $N>1$.  Since  $\cJ_\theta$ is unbounded from above, there is  $N_1>N+1$ such that the set
$\{N\le {\cJ}_\theta\le N_1,\ \inf_{t\le 1}e^{-V_t}\ge 1/N_1\}
$   
 is non-null. Then  
$$\Gamma^ \e:= \big\{N\le {\cJ}_\theta\le N_1,\  \inf_{t\le 1}e^{-V_t}\ge 1/N_1,\  |\zeta^\e_1|\le 1\big\}
$$ 
is also a non-null set for all sufficiently small $\e>0$. 

The process $P$ is not a subordinator and, therefore, we have only three possible cases.  

1) $\Pi_P(]-\infty,0[)>0$. Then  $\Pi_P(]-\infty,-\e_0[)>0$ for some $\e_0>0$.  
Due to  independence, the intersection of $\Gamma^\e$ with the set    
$$ 
\{| I_{\{x< -\e\}} x*\mu^P_1|\ge N_1(a_P^+N_1+N),\ I_{\{x >\e \}}*\mu_1^P=0\} 
$$
is non-null when $\e\in]0,\e_0[$. 
On this intersection we have that 
$$
Q_\theta\ge -a_P \cJ_\theta - \zeta^\e_1 - e^{-\theta V_{-}}\,I_{\{x< -\e\}} x*\mu^P_1\\
\ge -a_P^+N_1-1+a_P^+N_1+N\ge  N-1.
$$

2)  $\Pi_P(]-\infty,0[)=0$, $\Pi_P(h)=\infty$. Diminishing in the need  
$\e$ to ensure the inequality 
$\Pi_P(xI_{x>\e})\ge N_1(a_P^+N_1+N)$, we have on  the non-null set $\Gamma^\e\cap \{I_{x>\e}*\mu^P_1=0\}$ that  
$$
Q_\theta=-a_P \cJ_\theta - \zeta^\e_1 + e^{-\theta V_{-}}\,I_{\{x> \e\}} *\nu^P_1
\ge -a_P^+N_1-1+a_P^+N_1+N\ge N-1.
$$

3) $\Pi_P(]-\infty,0[)=0$, $\Pi_P(h)<\infty$, and 
$\Pi_P(h)-a_P>0$. Then on the non-null set $\{\cJ_\theta\ge N\}\cap \{I_{\{x>0\}}*\mu^P_1=0\}$ we have that 
$$
\cQ_\theta=
(\Pi_P(h)-a_P)\cJ_\theta\ge (\Pi_P(h)-a_P)N. 
$$

Since $N$ is arbitrary, in all three cases $Q_\theta$ is unbounded from above. \fdem

\begin{rem}
\label{II}
{\rm
If   $\cJ_1I_{\{V_1<0\}}$ is unbounded from above, so is  $Q_1I_{\{V_1<0\}}$.  
}
\end{rem}

\begin{rem} 
\label{II+} 
{\rm The proof above shows that in the case where $\sigma_P=0$  there is a constant $\kappa>0$ such that 
if the set $\{\cJ_\theta>N\}$ is non-null, then  $Q_\theta>\kappa N$ on its $\cF^{R,P}_1$-measurable non-null subset. 
The statement remains valid with obvious changes if the integration over the interval 
$[0,1]$ is replaced by the integral over arbitrary finite interval $[0,T]$. 
}
\end{rem}

\begin{lemm} \label{suppQ.1} 
$(i)$ The random variable  $\cJ_1$  is unbounded from above if and only if  $\sigma^2+\Pi(]-1,0[)>0$ or $\Pi(xI_{\{0<x\le 1\}})=\infty$. 

$(ii)$  The random variable  $\cJ_{-1}$ is unbounded from above if and only if  $\sigma^2+\Pi(]0,\infty[)>0$ or $\Pi(xI_{\{x<0\}})=-\infty$.
 \end{lemm} 
\noindent
{\sl Proof.}  In the  case where $\sigma^2>0$ the ``if" parts of the statements are obvious:  $W$ is independent of the  jump part of $V$ and  the distribution of the random variable  $\int_0^1e^{-\sigma \theta W_v}g(v)dv$, where $g>0$ is a deterministic function, has a support unbounded from above.  

So, suppose that $\sigma^2=0$ and consider  the ``if" parts  separately. 

\smallskip
$(i)$  
Let  $\Pi(]-1,0[)>0$, i.e. $\Pi(]-1,-\e[)>0$ for some $\e\in ]0,1[$.  Then the process 
$V$ given by (\ref{V}) admits the decomposition 
$$
V_t=at +  h*(\mu-\nu)_t+(\ln (1+x)-h)*\mu_t=(a-\Pi(xI_{\{-1<x\le -\e\}}))t+V'_t+V''_t, 
$$
where
\bean
V'_t&:=&
I_{\{ -\e< x\le 1\}}x*(\mu-\nu)_t+(\ln (1+x)-x)I_{\{-\e<x\le 1\}}*\mu_t\\
&&+\ln (1+x)I_{\{x> 1\}}*\mu_t,\\
V''_t&:=&\ln (1+x)I_{\{-1< x\le  -\e\}}*\mu_t. 
\eean
The processes $V'$ and $V''$ are independent. The decreasing process  $V''$ has jumps of the size not less than $|\ln (1-\e)|$ and the number of jumps on the interval $[0,t]$ is a Poisson random variable with parameter $t\Pi(]-1,-\e[)>0$.  Hence,  
 $V_t''$ is unbounded from below for any $t\in ]0,1[$. In particular, for any $N>0$, the set where  $e^{- V''}\ge N$ on the interval $[1/2,1]$  is non-null. 
The required property follows from these considerations. 

Let $\Pi(h(x)I_{\{x>0\}})=\infty$. We may assume without loss of generality that $\Pi(]-1,0[)=0$.  In this case, the process $V$  has only positive jumps.  Take arbitrary $N>1$ and choose 
$\e>0$ such that $\Pi(xI_{\{\e<x\le 1 \}})>2N$ and $\Pi(I_{\{0<x\le \e \}}\ln^2(1+x))\le 1/(32N^2)$.  We have
the decomposition 
$$
V_t=ct +V^{(1)}_t+V^{(2)}_t+V^{(3)}_t, 
$$
where the processes 
\bean
V^{(1)}&:=&I_{\{0<x\le \e\}}\ln (1+x)*(\mu-\nu),\\
 V^{(2)}&:=&I_{\{\e<x\le 1\}}\ln (1+x)*(\mu-\nu),\\
  V^{(3)}&:=& I_{\{x> 1\}}\ln (1+x)*\mu
\eean
are independent and the constant  
$$
c:= a+\Pi((\ln (1+x)-x)I_{\{0<x\le 1\}})<\infty. 
$$
By the Doob inequality  $P(\sup_{t\le 1}V^{(1)}_t < N/2)>1/2$. The processes  
$ V^{(2)}$ and $ V^{(3)}$ have no jumps on $[0,1]$ on  a non-null set. 
In the absence of jumps the trajectory of  $ V^{(2)}$ is the linear function 
$$
y_t=-\Pi(xI_{\{\e<x\le 1 \}})t\le -2Nt.
$$ 
It follows that $\sup_{1/2\le t\le 1}V_t\le c-N/2$
on the set of  positive probability.   This implies that  $\cJ_1$  
 is unbounded from above. 

\smallskip
$(ii)$  
Let  $\Pi(]0,\infty[)>0$, i.e. $\Pi(]0,\e[)>0$ for some $\e\in ]0,1[$.  Then 
$$
V_t=at +  h*(\mu-\nu)_t+(\ln (1+x)-h)*\mu_t=(a-\Pi(hI_{\{x>\e\}}))t+\tilde V'_t+\tilde V''_t, 
$$
where
\bean
\tilde V'_t&:=&
I_{\{ x\le \e\}}h*(\mu-\nu)_t+(\ln (1+x)-h)I_{\{x\le \e\}}*\mu_t,\\
\tilde V''_t&:=&\ln (1+x)I_{\{x>  \e\}}*\mu_t. 
\eean
The processes $\tilde V'$ and $\tilde V''$ are independent. The increasing process  $\tilde V''$ has jumps of the size not less than $\ln (1+\e)$ and the number of jumps on the interval $[0,t]$ is a Poisson random variable with parameter $t\Pi(]\e,\infty[)>0$.  Hence,  
 $V_t''$ is unbounded from above for any $t\in ]0,1[$.  In particular, for any $N>0$, the set where  $e^{V''}\ge N$ on the interval $[1/2,1]$ is non-null. These facts imply the required property. 

Let $\Pi(xI_{\{x<0\}})=-\infty$. 
We may assume without loss of generality that $\Pi(]0,\infty[)=0$. 
 In this case, the process $V$  has only negative  jumps.  Take arbitrary $N>1$ and choose 
$\e\in ]0,1/2[$ such that
 $$
 -\Pi(\ln(1+x)I_{\{-1/2<x\le -\e \}})>2N,
 \quad
 \Pi(I_{\{-\e < x<0 \}}\ln^2(1+x))\le 1/(32N^2).
 $$ 
This time we use the representation  
$$
V_t=\tilde ct +\tilde V^{(1)}_t+\tilde V^{(2)}_t+\tilde V^{(3)}_t, 
$$
where the processes
\bean
\tilde V^{(1)}&:=&I_{\{-\e<x<0\}}\ln (1+x)*(\mu-\nu),\\
\tilde V^{(2)}&:=&I_{\{-1/2<x\le-\e\}}\ln (1+x)*(\mu-\nu), \\
\tilde V^{(3)}&:=&I_{\{-1<x\le-1/2\}}\ln (1+x)*\mu
\eean
are independent and the constant 
$$
\tilde c:=a+\Pi(\ln (1+x)\,I_{\{-1/2<x<0\}}-h). 
$$
 By the Doob inequality  $\P(\sup_{t\le 1}\tilde V^{(1)}_t < N/2)>1/2$. The processes
$ \tilde V^{(2)}$ and $\tilde  V^{(3)}$ have no jumps on $[0,1]$  with strictly positive probability.  
In the absence of jumps the trajectory of  $ \tilde V^{(2)}$ is the linear function 
$$
y=-\Pi(\ln(1+x)I_{\{-1/2<x\le -\e \}})t\ge 2Nt. 
$$ 
It follows that $\sup_{1/2\le t\le 1}V_t\le \tilde c+N/2$
on a non-null set.   This implies that  $J_{-1}$   is unbounded from above. 

\smallskip
The ``only if'' parts of the lemma are obvious. \fdem 
\medskip

Summarizing, we conclude that $Q_1$ and $Q_{-1}$ (hence, $Y_{\infty}$) are unbounded from above 
if $\sigma^2>0$, or $\sigma^2_P>0$, or $P(|h|)=\infty$, or $\Pi(]-1,0[)>0$ and $\Pi(]0,\infty[)>0$. The remaining cases are treated in the following:    

\bigskip
\begin{lemm} 
\label{twocases}
Suppose that $\sigma=0$, $\Pi(|h|)<\infty$,  $\sigma_P=0$.  If $\Pi(]-1,0[)=0$ or $\Pi(]0,\infty[)=0$, then the random variable $Y_\infty$ is unbounded from above. 
\end{lemm}
\noindent
{\sl Proof.}  
By our  assumptions   $V_t=ct+L_t$ where  the constant $c:=a-\Pi(h)$, $\Pi\neq 0$,  and $L_t:=\ln (1+x)*\mu_t$.  
The assumption $\beta>0$ implies that $V_1<0$ with strictly positive probability and $V$ cannot be increasing or decreasing process. So, there are  two cases which we consider separately. 

$(i)$ 
$c<0$ and  $\Pi(]0,\infty[)>0$.  Take any $T>1$.  Then  the integral 
$\int_{[0,T]}e^{-V_{t-}}dt\ge T/e$ on the non-null set $\{L_T\le 1\}$.  By virtue of Remark  \ref{II+} on a non-null $\cF^{R,P}_T$-measurable  subset $\Gamma_T\subseteq \{L_T\le 1\}$ we have the bound  
$-\int_{[0,T]}e^{-V_{t-}}dP_t\ge K_T
$   
where $K_T\to \infty$ as $T\to \infty$. 
For every $T>1$  
$$
\P(\Gamma_T\cap 
\{ L_{T+1}-L_{T}\ge |c|(T+1)\})=\P(\Gamma_T) \P(L_{T+1}-L_{T}\ge |c|(T+1))>0. 
$$ 
Let  $\zeta^\e$ is the square integrable martingale  $\zeta^\e$  defined by  (\ref{zeta}) with $\theta=1$. 
Take $N>1$ sufficiently large and $\e>0$ sufficiently small to ensure that the set $\Gamma_T^{\e,N}$ defined as the intersection of  sets
$\Gamma_T\cap 
\{ L_{T+1}-L_{T}\ge |c|(T+1)\}$,  
$\big\{ \sup_{s\in [T,T+1]}e^{-V_s}\le N,\  \inf_{s\in [T,T+1]}e^{-V_s}\ge 1/N\big \}$,   and 
$\{ |\zeta^\e_{T+1}-\zeta^\e_{T}|\le 1\}$
is non-null.

 

\smallskip
Let us consider the representation 
\bean
Y_\infty&=&-\int_{[0,T]}e^{-V_{t-}}dP_t+a_P^\e\int_{]T,T+1]}e^{-V_{t-}} dt 
- \zeta_{T+1}^\e+\zeta_T^\e\\
&&-I_{]T,\infty[} e^{-V_{-}}xI_{\{|x|>\e\}}*\mu ^P_{T+1}
+ e^{-V_{T+1}}Y_{T+1,\infty}. 
\eean

Take arbitrary $y<0$ such that the set $\{Y_{T+1,\infty}>y\}$ is non-null.  

Since the process $P$ is not a subordinator with $\sigma_P=0$, it must satisfy one of the characterizing conditions 1), 2), 3) of Section \ref{prel}. Let us consider them consecutively.   

Suppose that  $\Pi_P(]-\infty,0[)>0$. Then there exists $\e_0>0$ such that $\Pi_P(]-\infty,-\e_0[)>0$. Due
to the independence, the intersection of $\Gamma_T^{\e,N}$ with the set    
$$ 
\tilde \Gamma_T^{\e,N}:=\{ I_{[T,\infty[} I_{\{x< -\e\}} *\mu^P_{T+1}\ge -(1/\e) N^2a_P^\e,\  I_{[T,\infty[}I_{\{x >\e \}}*\mu_{T+1}^P=0\} 
$$
is non-null when $\e\in]0,\e_0[$.

Due to independence, the intersection of   $\Gamma_T^{\e,N}\cap \tilde \Gamma_T^{\e,N}$ and $\{Y_{T+1,\infty}>y\}$
also is a non-null set. But on this intersection we have  inequality $Y_{\infty}\ge K_T-1+y$ implying that  
$Y_{\infty}$ is unbounded from above. 

Suppose that  $\Pi_P(]-\infty,0[)=0$, $\Pi_P(h)=\infty$. Thus,  for sufficiently small $\e>0$ we have $a_P^\e>0$. On the non-null
set 
$$
\Gamma_T^{\e,N}\cap \{I_{[T,\infty[}I_{\{x >\e \}}*\mu_{T+1}^P=0\} \cap  \{Y_{T+1,\infty}>y\}
$$
the inequality $Y_{\infty}\ge K_T-1+y$ holds  and we conclude as above. 

Finally, suppose that $\Pi_P(]-\infty,0[)=0$, $\Pi_P(h)<\infty$, and 
$\Pi_P(h)-a_P>0$. In this case we can use the representation 
\bean
Y_\infty&=&-\int_{[0,T]}e^{-V_{t-}}dP_t+(\Pi_P(h)-a_P)\int_{]T,T+1]}e^{-V_{t-}} dt 
\\
&&-I_{]T,\infty[} e^{-V_{-}}xI_{\{x>0\}}*\mu ^P_{T+1}
+ e^{-V_{T+1}}Y_{T+1,\infty}. 
\eean
On the non-null set $\Gamma_T^{\e,N}\cap \{I_{]T,\infty[} I_{\{x>0\}}*\mu ^P_{T+1}=0\}\cap \{Y_{T+1,\infty}>y\}$ we have that  
$Y_{\infty}\ge K_T+y$ implying that $Y_{\infty}$ is unbounded from above. 

\smallskip

$(ii)$ $c>0$ and $\Pi(]-1,0[)>0$. In this case there are  $\gamma,
\gamma_1\in ]0,1[$, $\gamma<\gamma_1$,  such that  
$\{I_{]-1,-\gamma[}*\mu_{1}= 0\}$, $\{I_{[-\gamma,-\gamma_1[}*\mu_{1/2}=I_{]-\gamma,-\gamma_1[}*\mu_{1}=N\}$,  and  $\{\ln(1+x)I_{]-\gamma_1,0[}*\mu_{1}\ge -1\}$ are non-null sets. Due to independence, their intersection $A_N$ is also non-null. 

On $A_N$ we have the bounds 
$$
c + N\ln (1-\gamma) -1\le V_1\le c + N\ln (1-\gamma_1)   
$$
and    
$$
\cJ_1:=\int_{[0,1]}e^{-V_{t-}}dt_t\ge e^{-c}\int_{[0,1/2]}e^{-\ln(1+x)*\mu_t}dt\ge \frac 12 e^{-c}(1-\gamma_1)^{-N}.  
$$ 
In virtue of Remark \ref{II+} there is a constant $\kappa_N$ an $\cF^{R,P}_1$-measurable non-null subset $B_N$ of $A_N$ such that 
$Q_1\ge \kappa_N$ on  $B_N$ and $\kappa_N\to  \infty$ as  $N\to \infty$.  

Take $T=T_N>0$ such that $cT+N\ln (1-\gamma) - 2\ge 0 $. 
The set $\{I_{]1,1+T[}\ln (1+x)*\mu_{1+T}\ge -1\}$ is non-null and 
its intersection  with $B_N$ is also non-null. On this intersection $e^{-V_{1+T}}\le 1$ and  
$$
c_1(N)\le V_{t-}\le c_2(N)   
$$
where $c_1(N):=c + N\ln (1-\gamma) -2$, $c_2(N):=c(T+1) + N\ln (1-\gamma_1)$. 


With this we accomplish the arguments by considering the cases 
corresponding to the properties 1), 2), and 3) with obvious modifications.  \fdem
\smallskip

With the above lemma the proof of Proposition \ref{YYY} is complete. \fdem 

\medskip
\noindent
{\bf Proof of the main theorem.} 
In view of (\ref{martbound}) and Lemma \ref{Q} we have that $\E\,|Q_1|^{\beta}<\infty$. 
The hypothesis on $\beta$ and Proposition \ref{YYY} allows us to use the results 
of the implicit renewal theory  on the tail behavior of distribution of $Y_\infty$ 
resumed in Theorem \ref{end} of Appendix. The reference to  Lemma  \ref{G-Paulsen} completes the proof.  \fdem 

\begin{rem} Note that the hypothesis
$\beta\in {\rm int}\,{\rm dom}\,H$ can be replaced by the slightly weaker assumption 
$\E e^{-\beta V_1}V_1^-<\infty$. 
\end{rem}

\begin{rem} The hypothesis $\cL(V_1)$ is non-arithmetic also can be replaced by a weaker one: one can assume that $\cL(V_T)$ is non-arithmetic for some $T>0$. Indeed,  due to the identity 
$\ln \E e^{-\beta V_T}=TH(\beta)$ the root $\beta$ does not depend on the choice of the time unit.   
\end{rem}

\small
The following lemma shows that the condition on $\cL(V_1)$ can be formulated in terms of the L\'evy triplets. 

\begin{lemm} The (non-degenerate) distribution of $V_1$ is arithmetic if and only if 
 $\sigma=0$,   $\Pi(\bbr)<\infty$, and there is $d>0$ such that $\Pi_V$ is concentrated on the lattice $\Pi (h)-a+\bbz d$.  
\end{lemm}
\noindent
{\sl Proof.} Recall   that  $\sigma_V=\sigma$   and $\Pi_V=\Pi\varphi^{-1}$ where $ \varphi: x\mapsto \ln (1+x)$. So, we have  $\Pi_V(\bbr)=\Pi(\bbr)$. If $\sigma_V>0$ or $\Pi_V(\bbr)=\infty$, the distribution of $V_1$ has a density,  see Prop. 3.12 in \cite{CT}. If $\sigma=0$  and $0<\Pi_V(\bbr)<\infty$,   then $V$ is a compound Poisson process with   drift $c=a-\Pi (h)$ and distribution of jumps $F_V:=\Pi_V/\Pi_V(\bbr)$. 
In such a case  $\cL(V_1)$ is concentrated on the lattice $\bbz d$  if and only if $\Pi_V$ is concentrated on the lattice $-c+\bbz d$. \fdem 

\section{Ruin with probability one} 
\label{prob1}

In this section we give conditions under which the 
 ruin is imminent whatever is the initial reserve.  

Recall the following  ergodic property of the  autoregressive 
process $(X_n^u)_{n\ge 1}$ with random coefficients (see, \cite{PeZe-06}, Prop. 7.1) which is defined   recursively by the relations 
\beq
\label{recur}
X_n^u=A_nX_{n-1}^u +B_n, \qquad n\ge 1, \quad X_0^u=u,  
\eeq
where $(A_n,B_n)_{n\ge 1}$ is a sequence of i.i.d. random variables in 
$\bbr^2$.   
\begin{lemm} 
Suppose that   $\E |A_n|^\delta<1$ and $\E |B_n|^\delta<\infty$ for some $\delta\in ]0,1[$. Then for any $u\in \bbr$ the sequence $X_n^u$ converges in $L^\delta$ (hence, in probability) to the random variable 
$$
X_\infty^0=\sum_{n=1}^\infty B_n\prod_{j=1}^{n-1}A_j
$$  
and for any bounded uniformly continuous function $f$
\beq\label{ergodic}
\frac 1N\sum_{n=1}^N f(X_n^u)\to \E f(X_\infty^0) \quad \hbox{in probability as } N\to \infty . \eeq  
\end{lemm}

Applying the lemma to the function $f(x)=I_{\{x<-1\}}-x I_{\{-1\le x<0\}}$ we get:  

\begin{coro} 
\label{coro}
Suppose that   $\E |A_n|^\delta<1$ and $\E |B_n|^\delta<\infty$ for some $\delta\in ]0,1[$. 

$(i)$ If $\P(X_\infty^0<0)>0$, then  $\inf_{n\ge 1}X_n^u<0$. 

$(ii)$ If $A_1>0$   and  $B_1/A_1$ is unbounded from below, then   $\inf_{n\ge 1}X_n^u<0$. 
\end{coro}

\noindent
{\sl Proof.} We get  $(i)$ by the straightforward application of (\ref{ergodic}) to the function $f(x):=I_{\{x<-1\}}-x I_{\{-1\le x<0\}}$.  The statement $(ii)$ follows from $(i)$.  Indeed, put  $X_\infty^{0,1}:=\sum_{n= 2}^\infty B_n\prod_{j= 2}^{n-1}A_j$. Then 
$$
X_\infty^0=B_1+A_1\, X_\infty^{0,1}=A_1(X_\infty^{0,1}+B_1/A_1). 
$$ 
Since $B_1/A_1$ and $X_\infty^{0,1}$ are independent and the random variable $B_1/A_1$ is unbounded from below, $\P(X_\infty^0<0)>0$.  \fdem

\smallskip
Let $M_j$ and $Q_j$ be the same as in  (\ref{MQ}). 

\begin{prop} 
\label{Psi=1}
 Suppose that   $\E M_1^{-\delta}<1$ and $\E M_1^{-\delta}| Q_1|^\delta<\infty$ for some $\delta\in ]0,1[$.  If $Q_1$ is unbounded from above, then $\Psi (u)\equiv 1$.  
\end{prop} 
\noindent{\sl Proof.} 
The process $X^u$ solving the equation (\ref{basic}) and restricted to the integer values of the time scale  admits  the representation 
$$
X_n^u=e^{V_n-V_{n-1}}X_{n-1}^u +e^{V_n} \int_{]n-1,n]}e^{-V_{t-}}dP_t, \qquad n\ge 1, \quad X_0^u=u.     
$$
That is,  $X_n^u$ is given by (\ref{recur}) with $A_n=M_n^{-1}$
and $B_n=-M_n^{-1}Q_n$. The result follows from the statement $(ii)$ of  Corollary \ref{coro}. 
\fdem

Now we give more specific conditions of the ruin with probability one  in terms of the triplets.  


\smallskip

%
%
%
%
%

\begin{theo} 
\label{ruintheo}
Suppose that  $0\in {\rm int}\,{\rm dom}\,H$ and $\Pi_P(|\bar h|^\e)<\infty$ for some $\e>0$. If $a_V+\Pi(\bar h(\ln (1+x)))\le 0$, then $\Psi (u)\equiv 1$. 
\end{theo} 
\noindent
{\sl Proof.}  Note that $D^-H(0)=-a_V-\Pi(\bar h(\ln (1+x)))$. If $D^-H(0)>0$, then 
 for all $q<0$ sufficiently close to zero $H(q)<0$, i.e. $\E M_1^{q}<1$. 
By virtue of Lemma \ref{Q/M} $\cL(M_1^{-1}Q_1)=\cL(Q_{-1})$. 
If $\Pi_P(|\bar h|^\e)<\infty$ for some $\e>0$, then the same arguments as in the proof of the first part of the proof of Lemma \ref{Q} lead to the conclusion that  $\E |Q_{-1}|^{q}<\infty$ for sufficiently small $q>0$.  
 To get the result we can use Proposition \ref{Psi=1}.  Indeed, 
by virtue of Lemmata \ref{suppQ} and \ref{suppQ.1}(i) the random variable 
 $Q_1$ is   unbounded from above except, eventually, the  case where 
 $\sigma^2=0$, $\sigma^2_P=0$, $\Pi(]-1,0[)=0$,  $\Pi(xI_{\{0<x\le 1\}})<\infty$, and $\Pi\neq 0$. But under such constrains on the characteristics 
 the distribution of $X_\infty^0$ coincides with the distribution of the integral $\int_0^\infty e^{V_s}dP_s$.  Using the arguments similar to those  in the proof of Lemma   \ref{twocases}(i), it is easy to prove that the latter charges $]-\infty,0[$ and 
  we can  apply  Corollary \ref{coro}$(i)$.
  
\smallskip


In the case where $D^-H(0)= 0$ we consider, following \cite{PeZe-06}, the  discrete-time process  $(\tilde X^u_n)_{n\in {\bbn}}$ where   $\tilde X^u_n=X_{T_n}$
and  the descending ladder times $T_n$ of the random walk $(V_n)_{n\in {\bbn}}$ which are defined as follows: $T_0:=0$, 
$$
T_n:=\inf\{k>T_{n-1}\colon\ V_k-V_{T_{n-1}}<0\}. 
$$
Since $J(q)=\Pi \big(I_{\{|\ln (1+x)|>1\}}(1+x)^{-q}\big)<\infty$ for any $q\in ]\underline q,\bar q[$, we have that $\Pi(\ln^2(1+x)))<\infty$. It follows that  the formula (\ref{V}) can be written as 
$$
V_t=\big(a-\sigma^2/2 -\Pi(h)\big)t + \sigma W_t+ \ln (1+x)*\mu_t,
$$
$\E V_1^2<\infty$, and the condition $D^-H(0)= 0$ means that $\E V_1=0$. 

Accordingly to Theorem 1a in Ch. XII.7 of Feller's book \cite{Feller1971} and the remark preceding the citing  theorem, 
the above properties imply that there is a finite constant $c$ such that 
\begin{equation}\label{sec:Ras.5}
\P\left(T_{1} > n\right) \le cn^{-1/2}.
\end{equation}
It follows, in particular, that the differences $T_n-T_{n-1}$ are well-defined and form  a sequence of finite  independent random variables distributed as $T_1$. The discrete-time process $\tilde X^u_n=X^u_{T_n}$ 
has the representation 
$$
\tilde X_n^u=e^{V_{T_n}-V_{T_{n-1}}}\tilde X_{n-1}^u +e^{V_{T_n}} \int_{]T_{n-1},T_n]}e^{-V_{t-}}dP_t, \qquad n\ge 1, \quad \tilde X_0^u=u,     
$$
and solves the linear equation 
$$
\tilde X_n^u=\tilde A_n\tilde X_{n-1}^u +\tilde B_n, \qquad n\ge 1, \quad X_0^u=u,  
$$
where
$$
\tilde A_n:=e^{V_{T_n}-V_{T_{n-1}}}, \qquad \tilde B_n:= e^{V_{T_n}} \int_{]T_{n-1},T_n]}e^{-V_{t-}}dP_t,  
$$
and $\tilde B_1/\tilde A_1=Y_{T_n}$ where $Y$ is given by (\ref{Y_t}). 

By construction, $ \tilde A^\delta_1<1$ for any $\delta>0$. 

Using the definition of $Q_j$ given by  (\ref{MQ}) we have that 
$$
\vert \tilde B_{1} \vert
\le  \sum^{T_{1}}_{j=1}\,e^{V_{T_{1}}-V_{j-1}}
\,\vert {Q}_{j}\vert\le \sum^{T_{1}}_{j=1}\vert Q_{j}\vert .
$$
According to Lemma \ref{Q}   $\E\vert Q_{1}\vert^{p}<\infty$ 
 for some $p\in ]0,1[$.
Then for  $r\in ]0,p/5[$ and $l_{n}:=[n^{4r}]$, we have, using the Chebyshev inequality and  (\ref{sec:Ras.5}), 
that 
\begin{align*}
\E\,\vert \tilde B_{1} \vert^{r}&\le 1+r\sum_{n\ge 1}\,{n^{r-1}}\P\left(  \sum^{T_{1}}_{j=1}
\,\vert Q_{j}\vert
>n
\right)\\[2mm]
&\le  1+r\sum_{n\ge 1}\,{n^{r-1}}\,\P\left(  \sum^{l_{n}}_{j=1}
\,\vert Q_{j}\vert>n\right)+ r\sum_{n\ge 1}\,{n^{r-1}}
\P\left(  T_{1}>l_{n}
\right)
\\[2mm]
&\le 
 1+r\E\vert Q_{1}\vert^{p}\sum_{n\ge 1}\,l_{n}n^{r-1-p} +rc\sum_{n\ge 1}\,n^{r-1}l_{n}^{-1/2}<\infty.
\end{align*}
 
To apply Corollary \ref{coro}$(ii)$ it remains to check that  $Y_{T_1}$ is unbounded from above.  Since  $\left\{Q_1>N\,,
\,V_{1}<0\right\}\subseteq \{Y_{T_1}>N\}$, it is sufficient to check that the probability of the set in the left-hand side is strictly positive  for all $N>0$, or,  by virtue of Remark  \ref{II},  that 
\beq
\label{IN}
\P(\cJ_1>N,\; V_1<0)>0 \qquad \forall\; N>0. 
\eeq

 Let $\sigma^2>0$.  Taking into account that the conditional distribution of the process $(W_s)_{s\le 1}$ given $W_1=x$ is the same as the (unconditional) distribution  of the Brownian bridge $B^x=(B_s^x)_{s\le 1}$ with  $B^x_s =W_s+s(x-W_1)$  we easily get that for any bounded  positive  function 
 $g$ and  any $y,M \in\bbr$ the  probability 
$$
\P\left( 
\int_0^1\,e^{-\sigma W_{v}}g(v)d v>y
\,,\,
W_{1}<M\right)>0,   
$$
cf. with  Lemma 4.2 in \cite{KP}. 
This implies  (\ref{IN}). 

Suppose that $\sigma^2=0$, but $\Pi(]-1,0[)>0$, i.e. $\Pi(]-1,-\e[)>0$ for some 
$\e\in ]0,1[$. In the decomposition $V=V^{(1)}+V^{(2)}$, where
\bean
V^{(1)}_t&=&I_{\{-1<x\le -\e\}}\ln (1+x)*\mu_t,\\
V^{(2)}_t&=&(a-\Pi(hI_{\{-1<x\le -\e\}}))t+I_{\{x> -\e\}}h*(\mu-\nu)_{t}\\
&&+I_{\{x> -\e\}}(\ln (1+x)-h)*\mu_t,
\eean
the processes $V^{(1)}$ and $V^{(2)}$ are independent.  The process $V^{(1)}$ is decreasing by negative jumps
whose absolute value are larger or equal than $|\ln (1-\e)|$ and the number of  jumps on the interval    $[0,1/2]$ has the  Poisson distribution  with parameter $(1/2)\Pi(]-1,-\e[)>0$. 
Thus,  $\P (V^{(1)}_{1/2}< -n)>0$ for any real $n$.  It follows that 
\bean
\P(\cJ_1>N,\, V_1<0)&\ge& \P\Big(\int_0^1e^{-V_t} dt>N,\, V_1<0,\, V^{(1)}_{1/2}< -n \Big)\\
&\ge&  \P\Big(e^n\int_{1/2}^{1}e^{-V^{(2)}_t} dt>N,\, V^{(2)}_1<n,\, V^{(1)}_{1/2}< -n \Big)\\
&=&  \P\Big(\int_{1/2}^{1}e^{-V^{(2)}_t} dt>Ne^{-n},\, V^{(2)}_1<n\Big)\P( V^{(1)}_{1/2}< -n ). 
\eean
The right-hand side is strictly positive for sufficiently large $n$ and   (\ref{IN}) holds.

The case where $\Pi(xI_{\{0<x\le 1\}})=\infty$ is treated similarly as in the last part of the proof of Lemma \ref{suppQ.1}$(i)$.   

The exceptional case  is treated by a reduction to Corollary \ref{coro}$(i)$. 
\fdem

\smallskip The above theorem implies that in the classical model with negative risk sums (where $\sigma_P=0$, the jumps of $P$ are positive and form a compound Poisson process, $\Pi_P(|x|)<\infty$, trend is negative, i.e. $a_P-\Pi_P(x)<0$) and investments into a risky asset with the price following a geometric Brownian motion (that is, $\Pi=0$ and  $\sigma\neq 0$), the ruin is imminent if $a_V= a-\sigma^2/2\le 0$.  

\section{Examples}
\label{Sec6}
\noindent
{\bf Example 1.}  Let  us consider the  model with negative risk sums in which 
$\Pi_P(dx)=\lambda F_P(dx)$ where the constant $\lambda>0$ and  the probability distribution $F_P(dx)$ is concentrated on $]0,\infty[$,  and 
$$
a_P^0:=\lambda \int_{[0,1]}x F_P(d x)-a_P. 
$$ 
The  process $P$ admits the representation as sum of an independent Wiener process  with drift and a compound Poisson process: 
 \begin{equation}
  \label{risk.1}
 P_{t}=-a_P^0t +\sigma_PW_t^P+\sum^{N^P_{t}}_{j=1}\,\xi_{j},
 \end{equation}
where the Poisson process $N^P$  with 
intensity $\lambda_P$ is independent of  the sequence 
$(\xi_j)_{j\ge 1}$ of  positive i.i.d. random variables with common
distribution $F_P$.   

Suppose that the price process   is a geometric Brownian motion
$$
\cE_t(R)=e^{V_t}=e^{(a-\sigma^2/2)t+\sigma W_t},
$$
that is,  $\sigma\neq 0$, $\Pi=0$. 

For this model $\underline q=-\infty$, $\bar q=\infty$. The condition $D^+H(0)<0$ is reduced to the inequality  $\sigma^2/2<a$ and the function  $H(q)=(\sigma^2/2-a+q\sigma^2/2)q$ has the root $\beta=2a/\sigma^2-1>0$. Suppose that $\sigma^2_P+(a_P^0)^+>0$.
 By Theorem \ref{Main2+} the exact asymptotic $\Psi(u)\sim C_{\infty}u^{-\beta}$, as $u\to \infty$, 
 holds if  $\E\xi_1^{\beta}<\infty$. Since the exponential distribution has the above property,  we recover, as a very particular case the  asymptotic result  of \cite{KP} where  it was assumed that $\sigma^2_P=0$ and $a_P^0>0$. 
 
\smallskip 
If $\sigma^2_P+(a_P^0)^+>0$,  $\sigma^2/2\ge a$,  and $\E\xi_1^{\epsilon}<\infty$ for some $\epsilon>0$, then Theorem 
\ref{ruintheo}  implies that $\Psi(u)\equiv 1$. 
 
\smallskip
\noindent
{\bf Example 2.}  Let the process $P$ be again  given by (\ref{risk.1}) 
 and suppose that  
the price  process  has a jump component, namely,  
$$
\cE_t(R)=\exp\Big\{(a-\sigma^2/2)t+\sigma W_t+\sum_{j=1}^{N_t}\ln (1+\eta_j)\Big\}, 
$$ 
where the Poisson process $N$  with 
intensity $\lambda>0$ is independent on the sequence 
$(\eta_j)_{j\ge 1}$ of  i.i.d. random variables with common
distribution $F$ not concentrated at zero and $F(]-\infty,-1])=0$, see \cite{LL}, Ch. 7.  That is, the log price process is represented as 
$$
V_t=(a-\sigma^2/2)t+\sigma W_t+ \ln (1+x)*\mu_t,  
$$
 where $\Pi(dx)=\lambda F(dx)$.  The function $H$ is given by the formula 
 $$
 H(q)=(\sigma^2/2-a+q\sigma^2/2)q+\lambda (\E\, (1+\eta_1)^{-q}-1).  
 $$
Suppose that $\E\, (1+\eta_1)^{-q}<\infty$ for all $q>0$. Then $\bar q=\infty$. 


Let  $\sigma\neq 0$.  Then  $\limsup_{q\to \infty} H(q)/q=\infty$.  If   
\beq
\label{DH0}
D^+H(0)=\sigma^2/2-a -\lambda  \E\ln (1+\eta_1)<0,
\eeq
then the root $\beta>0$ of the equation $H(q)=0$ does exist. Thus, if   $\E\xi_1^{\beta}<\infty$, then Theorem \ref{Main2+} can be applied to get 
that  $\Psi(u)\sim C_{\infty}u^{-\beta}$ where $C_{\infty}>0$. 

If $\E(1+\eta_1)^{-\beta_1}<1$ (resp., $\E(1+\eta_1)^{-\beta_1}>1$), the root $\beta$ is smaller (resp., larger) than  $2a/\sigma^2-1$, the value of the root of $H$ in model of the first example where the price process is continuous. 

\smallskip
Let  $\sigma= 0$.  If 
  $$
D^+H(0)=-a -\lambda \E\, \ln (1+\eta_1)<0,
$$
and 
$$\limsup_{q\to \infty} q^{-1}\E\, \big((1+\eta_1)^{-q}-1\big)>a/\lambda,
$$ 
then the root $\beta>0$ also exists.  Theorem \ref{Main2+} can be applied  
when   $0<\P(\eta_1>0)<1$ and the we have  exact asymptotic  
if the distribution of $\ln (1+\eta_1)$ is non-arithmetic.  

\smallskip
Suppose that $\E\, (1+\eta_1)^{-q}<\infty$ for all $q\in \bbr$. Then $\underline q=-\infty$, $\bar q=\infty$. 
If $\sigma^2/2-a -\lambda  \E\ln (1+\eta_1)\ge 0$,  $\sigma^2+\P(\eta_1<0)>0$, and   $\E |\xi_1|^\e<\infty$ for some $\e>0$,  then  $\Psi(u)\equiv 1$ in virtue of Theorems \ref{ruintheo}.

\section{Appendix: tails of distributions solving distributional equations}
\label{sec:A}
\setcounter{equation}{0}
\label{App}

\subsection{Kesten--Goldie theorem}

Here we present a short account of needed results on   distributional equations (random equations in the terminology of \cite{Go91})  
\begin{equation}\label{3.1}
Y_{\infty} \stackrel{d}{=}Q+M\,Y_{\infty},\quad
Y_{\infty} \ \ \mbox{independent of}\ (M,Q),
\end{equation}
where $(M,Q)$ is a given two-dimensional random variable with $M>0$ and $\P(M\neq 1)>0$ and   $\stackrel{d}{=}$ is the equality in law.  This is a symbolical notation  which means that we are given in fact  a two-dimensional 
distribution $\cL$ on $\bbr\times ]0,\infty[$ not concentrated 
on $\bbr\times \{1\}$
 and the problem is to find a probability space with random variables $Y_{\infty} $ and $(M,Q)$ on it such that $Y_\infty$ and $(M,Q)$ are independent,  
$\cL(M,Q)=\cL$,  
 and  $\cL(Y_{\infty})= \cL(Q+M\,Y_{\infty})$.
The uniqueness  in this problem means the uniqueness of the distribution of $Y_{\infty} $.   

\smallskip
In the sequel  $(M_j,Q_j)$ will be an i.i.d. sequence whose generic 
term $(M,Q)$ has the distribution $\cL$ and $Z_j:=M_1\dots M_j$, $Z_n^*:=\sup_{j\le n}Z_j$. 
\smallskip 
 
If there is $p>0$ such that $\E M^p<1$ and $\E |Q|^p<\infty$, then the  solution 
$Y_\infty$ of (\ref{3.1})  can be easily realized on the probability space $(\Omega,\cF,\P)$ where the sequence    
$(M_j,Q_j)$ is defined --- just as the limit in $L^p$ of the series $\sum_{j\ge 1}Z_{j-1}Q_j$,  see the beginning of the proof of Proposition \ref{YYY}.  

\smallskip

The following classical result of the renewal theory is the Kesten--Goldie theorem, see Th. 4.1 in \cite{Go91}: 
\begin{theo} 
\label{Le.3.2}
Suppose that $(Q,M)$ is such that the distribution of   $\ln\,M$ 
is non-arithmetic and, for  some $\beta>0$,
\begin{align}\label{3.3}
\E\,M^\beta=1, \ \ \ \E\,M^\beta\,(\ln\,M)^+<\infty, \  \  \ \E\,|Q|^\beta<\infty. 
\end{align}
Then 
\bean
&&\lim_{u\to\infty}\,u^\beta\,\P(Y_{\infty} >u)=C_+<\infty,\\
&&\lim_{u\to\infty}\,u^\beta\,\P(Y_{\infty} <-u)=C_-<\infty,
\eean
where $C_++C_->0$. 
\end{theo} 

\smallskip    
Theorem \ref{Le.3.2} left open the question when the constant $C_+$ is strictly positive. 
Recently,  
 Guivarc'h and  Le Page showed for the above case where the distribution of $\ln M$ is non-arithmetic that $C_+>0$ if and only if $Y_{\infty} $ is unbounded from above, 
 see \cite{GL}  and also the  paper \cite{BD} for simpler  arguments. The remaining part of the appendix deals mainly with the arithmetic case. 

\subsection{Grincevi\u{c}ius theorem}

The theorem below is a simplified  version of  Th.2(b), \cite{Grincevicius1975_a}, but with a slightly weaker assumption  on $Q$, 
namely, $\E|Q|^\beta<\infty$, used in our study. For the reader convenience we give its complete proof after recalling  some  concepts and facts from the renewal theory. 

\begin{theo} 
\label{Le.3.2++}
Suppose that (\ref{3.3}) holds 
 and the distribution of $\ln M$ is concentrated on the lattice $\bbz d$ where $d>0$. Then 
\beq
\label{3.3++4+}
\limsup_{u\to\infty}\,u^\beta\,\P(Y_{\infty} >u)<\infty.
\eeq
\end{theo}

\smallskip

We consider the convolution-type linear operator which is well-defined for all positive as well as for  (the Lebesgue)  integrable functions
by  the formula 
\begin{equation}
\label{transform++}
\check{\psi}(x)= \int^{x}_{-\infty}\,e^{-(x-y)}\,\psi(y) d y.
\end{equation}
Clearly,  the functions $\psi$ and $\check \psi$  are integrable or not 
simultaneously and  
$$
\int_{\bbr}\check{\psi}(x)d x= 
\int_{\bbr}\psi(x)d x.
$$
Suppose that  $\psi\ge 0$  is  integrable. Then   $\check{\psi}(x+\delta)\ge e^{-\delta}\check{\psi}(x)$  for any $\delta>0$ and 
$$
\delta \inf_{x\in [j\delta , (j+1)\delta]} 
\check{\psi}(x)\ge  \delta e^{-\delta} \check{\psi}(j\delta)\ge e^{-2\delta} \int_{(j-1)\delta}^{j\delta}\check{\psi}(x)dx
$$
implying that   
$$
\underline U(\check{\psi},\delta):=\delta\sum_{j\in {\mathbb Z}}\inf_{x\in [j\delta , (j+1)\delta]} 
\check{\psi}(x)\ge  e^{-2\delta}\int_{\bbr}\check{\psi}(x)dx.
$$
Similarly, 
$$
\bar U(\check{\psi},\delta):=\delta\sum_{j\in {\mathbb Z}}\sup_{x\in [j\delta , (j+1)\delta]} 
\check{\psi}(x)\le e^{2\delta}\int_{\bbr}\check{\psi}(x)dx.
$$
Thus, $\bar U(\check{\psi},\delta)<\infty$ and   $\bar U(\check{\psi},\delta)- \underline U(\check{\psi},\delta)  \to 0$ as $\delta\to \infty$. These two properties mean, by definition, that  
the function $\check{\psi}$ is directly Riemann integrable. Arguing with the positive and negative parts, we obtain that  if $\psi$ is integrable, then 
$\check{\psi}$ is directly Riemann integrable. 

We shall use in the sequel the following renewal theorem for the random walk $S_n:=\sum_{i=1}^n\xi_i$ on a lattice, see Prop.  2.1, \cite{IksanovPolotskiy2016}.   
\begin{prop} 
\label{Iks}
Let   $\xi_i$ be  i.i.d. random variables  taking values  in the lattice $\bbz d$, $d>0$,   and having finite expectation $m:=\E\xi_i>0$. Let 
 $F:\bbr\to \bbr$ be a measurable  function. If $x\in \bbr$ is such that  $\sum_{j\in \bbz}|F(x+jd)|<\infty$, then 
$$
\lim_{n\to \infty}\E \sum_{k\ge 0}F(x+nd-S_k)=\frac dm \sum_{j\in \bbz}F(x+jd).
$$

\end{prop}

\smallskip
\noindent{\sl Proof of Theorem \ref{Le.3.2++}.}
Let the solution of (\ref{3.1}) be realized on some probability space $(\Omega,\cF,\P)$. We shall use the notation $(M,Q)$ instead of  $(M_1,Q_1)$. 
Put $ \bar G(u):=\P(Y_{\infty}> u)$ and  $g(x):=e^{\beta x} \bar G(e^x)$. Since $Y_\infty$ and $M$ are independent, $\P(MY_{\infty}>e^{x})=\E \bar G(e^{x-\ln M})$. Defining the new probability measure $\tilde \P:=M^{\beta} \P$ and noting that 
$$
e^{\beta x}\P(MY_{\infty}>e^{x})
=\E M^\beta e^{\beta(x-\ln M)} \bar G(e^{x-\ln M})= \tilde \E g(x-\ln M)
$$
we obtain the following identity (called {\it renewal equation}): 
\begin{equation}
\label{ren.eq_+++}
g(x)=D(x)+\tilde \E g(x-\ln M),
\end{equation}
where $D(x):=e^{\beta x}\left(\P(Y_{\infty} >e^{x})-\P(MY_{\infty}>e^{x})\right)$. 
The  Jensen 
inequality for the convex function $x\mapsto x\ln x $ implies that 
$\tilde \E \ln M=\E M^\beta \ln M>0$ 
and, hence, $\tilde \E |\ln M|<\infty$.  

Let us check that the function $x\mapsto D(x)$ is integrable. 
To this aim, we note that for any random variables $\xi,\eta$
$$
|\P(\xi >s)-\P(\eta>s)|\le \P(\eta^+\le s<\xi^+)+\P(\xi^+\le s<\eta^+).  
$$
Using the Fubini theorem we obtain that  
$$
 \int_0^\infty \P(\eta^+\le s<\xi^+)s^{\beta-1}ds=  \E I_{\{\eta_+<\xi_+\}}\int_{\eta^+}^{\xi^+}s^{\beta-1}ds=\frac 1\beta \E\big((\xi^+)^\beta-(\eta^+)^\beta)\big)^+. 
$$
Applying this bound with $\xi:=Q+MY_\infty\stackrel{d}{=}Y_{\infty}$ and $\eta:=M Y_{\infty}$ we get that 
$$
\int_{\bbr}\,\vert D(x)\vert d x=\int_0^\infty|\P(\xi >s)-\P(\eta>s)|s^{\beta-1}ds\le \frac 1\beta \E\big\vert(\xi^+)^\beta-(\eta^+)^\beta)\big\vert  
$$
and it remains to verify  that  
\beq
\label{bbound}
\E|((Q+\eta)^+)^\beta-(\eta^+)^\beta|
<\infty
\eeq
when  $\E \vert Q\vert^{\beta}<\infty$. 
But  
$|((Q+\eta)^+)^\beta-(\eta^+)^\beta|=\zeta_1+\zeta_2$ with positive  summands 
\bean
\zeta_1&:=&I_{\{-Q<\eta\le 0\}}(Q+\eta)^\beta+I_{\{0<\eta\le -Q\}}\eta^\beta \le |Q|^\beta, \\
\zeta_2&:=&I_{\{Q+\eta>0,\,\eta >0\}}|(Q+\eta)^\beta-\eta^\beta|.
\eean
If $\beta\le 1$, then $\zeta_2$ is also dominated by $|Q|^\beta$. If 
$\beta>1$, then  the  inequality $|x^\beta-y^\beta|\le 
\beta |x-y|(x\vee y)^{\beta-1}$ for $x,y\ge 0$ combined with  
the inequality $(|a|+|b|)^{\beta-1}\le 2^{(\beta-2)^+}(|a|^{\beta-1}+|b|^{\beta-1})$ leads to the estimate 
$$
\zeta_2\le 2^{(\beta-2)^+}\beta |Q| (|\eta|^{\beta-1}+|Q|^{\beta-1}). 
$$
Using the independence of $(M,Q)$  and $Y_\infty$, the  H\"older inequality,  and taking into account that $\E M^\beta=1$ and   $\E| Y_\infty|^p<\infty$ for $p\in [0,\beta[$ we get that   
$$
\E|Q||\eta|^{\beta-1}=\E|Q|M^{\beta-1}\E| Y_\infty|^{\beta-1}\le (\E|Q|^\beta)^{1/\beta}\E |Y_\infty|^{\beta-1} <\infty. 
$$
Thus, (\ref{bbound}) holds. 
     
     The integrability of $D$ allows us to transform (\ref{ren.eq_+++}) into the 
     equality 
     $$
\check{g}(x)=\check{D}(x)+\tilde \E \check{g}(x-\ln M). 
$$
Iterating it, we obtain that  
\beq
\label{hatv}
\check{g}(x)=\sum_{n= 0}^{N-1}\tilde \E\check{D}(x-{S}_{n})+\tilde \E\check{g}(x-{S}_{N}), 
\eeq
%
%
%
where $S_{0}=0$ and  $S_{n}:=\sum_{i=1}^n\xi_i$ for $n\ge1$,  $(\xi_{i})$
is a sequence of independent random variables on $(\Omega,\cF,\tilde \P)$  independent on $Y_{\infty}$ such that the distribution  $\cL(\xi_i,\tilde \P)=\cL(\ln M,\tilde \P)$. In particular, $\tilde \E e^{-\beta \xi_i}=1$. 

By the strong law of large numbers $S_N/N\to \tilde \E \ln M>0$  $\tilde \P$-a.s., $N\to \infty$, and, therefore,  $y- S_N\to -\infty$ $\tilde \P$-a.s. for every $y$.   
Since $\tilde \E e^{-\beta S_N}=1$, we have by dominated convergence that  
$$
\tilde \E g(y-{S}_{N})=\tilde \E e^{\beta (y-S_N)}\bar G(e^{y- S_N})\to 0.   
$$
It follows that the remainder term $\tilde \E\check{g}(x-{S}_{N})$ in (\ref{hatv}) tends to zero,  thus,   
\beq
\label{hatv+}
\check{g}(x)=\sum_{k\ge  0}\tilde \E\check{D}(x-{S}_{k}). 
\eeq

Using 
 Proposition \ref{Iks} (with $F=\check D$) we obtain that for any $x>0$
 \begin{equation}
\label{limit_ren_th}
 \lim_{n\to\infty}\check{g}(x+d n)=
 \frac{d}{\tilde \E \ln M}\sum_{j\in {\mathbb Z}}
 \check{D}(x+jd)
 \le 
 \bar U(\check{D},d)
 <\infty.
\end{equation}
Replacing in the  integrant the function $\bar G(e^y)$ by its smallest value $G(e^{x})$ we obtain that  
$$
\check{g}(x):=\int^{x}_{-\infty}\,e^{-(x-y)}e^{\beta y} \bar G(e^y) d y
\ge \frac{1}{\beta+1}g(x)
$$
and, therefore, 
$$
\limsup_{u\to\infty}\,u^{\beta}\P(Y_{\infty}> u)=
\limsup_{x\to\infty}\,g(x)
\le (\beta+1)
\limsup_{x\to\infty}\,\check{g}(x)<\infty.
$$
Theorem \ref{Le.3.2++} is proven.
\fdem



\subsection{Buraczewski--Damek approach}

The following result, usually formulated in terms of the supremum 
of the random walk $S_n:=\sum_{i=1}^n{\ln M_i}$, 
 is well-known (see, e.g., Th. A, \cite{Kesten} for much more general setting).  

\begin{prop}
\label{Pr.prod_M_j}  
If $M$ satisfies (\ref{3.3}), then  
\beq 
 \label{prod_r_vs}
 \liminf_{u\to\infty}u^{\beta}
 \P (Z^*_{\infty}>u)>0.
 \eeq
\end{prop} 
\noindent{\sl Proof.} Let $F(x):= \P(\ln M\le x)$, $\bar F(x):=1-F(x)$,  
$S_n:=\sum_{i=1}^n{\xi_i}$ where $\xi_i:=\ln M_i$.  The function  $\bar H(x):=
\P(\sup_n S_{n}>x)$ admits the representation  
$$
\bar H(x)=\P(\xi_{1}>x)+\E\,I_\zs{\{\xi_{1}\le x\}}\,
\bar H(x-\xi_{1})
=\bar F(x)
+
\,
\int^{x}_{-\infty}\,
\bar H(x-t)d F(t).
$$
Putting   $Z(x):=e^{\beta x}\bar H(x)$, $z(x):=e^{\beta x}\bar F(x)$, and $\tilde \P:=e^{\beta\xi_1} \P$,  we obtain from here that
\beq\label{ren_eq}
Z(x)=z(x)+
\tilde \E\,Z(x-\xi_1)I_{\{\xi_1\le x\}}.
\eeq
The same arguments as were used in deriving (\ref{hatv}) lead
to the representation 
\begin{equation}
\label{ren_fun}
Z(x)=\tilde \E\,\sum_{k\ge 0}\,z(x-S_{k})I_{\{S_k\le x\}}.
\end{equation}
The function $\hat z(x):=z(x)I_{\{x\ge 0\}}$ is directly Riemann integrable. Indeed, for $j\ge 0$ we have that 
$$
\sup_{x\in [j\delta , (j+1)\delta]} 
z(x)\le e^{\beta (j+1)\delta}\,
\bar{F}(j\delta )\le e^{2\beta \delta}\int^{j\delta}_{(j-1)\delta}
e^{\beta v}
\bar{F}(v)d v
$$
and, therefore, 
$$
\bar U(\hat z,\delta)=\delta z(0)+\delta\sum_{j\ge 0}\sup_{x\in [j\delta , (j+1)\delta]} 
z(x)\le \delta z(0)+ e^{2\beta \delta}\int_{-\delta}^\infty e^{\beta v}
\bar{F}(v)d v.
$$
In the same spirit 
$$
\inf_{x\in [j\delta , (j+1)\delta]} 
z(x)\ge e^{\beta j\delta}\,
\bar{F}((j+1)\delta )\ge e^{-2\beta \delta}\int^{(j+2)\delta}_{(j+1)\delta}
e^{\beta v}
\bar{F}(v)d v
$$
and 
$$
\underline  U(\hat z,\delta)=\delta\sum_{j\ge 0}\sup_{x\in [j\delta , (j+1)\delta]} 
z(x)\ge e^{-2\beta \delta}\int_{\delta}^\infty e^{\beta v}
\bar{F}(v)d v.
$$

Taking into account that 
$$
\int_{\bbr}e^{\beta v}
\bar{F}(v)d v
= \frac{1}{\beta}\,\E\,e^{\beta\xi_{1}}
=\frac{1}{\beta}<\infty.
$$
We get from here that $\bar U(\hat z,\delta)<\infty$   and 
$\bar U(\hat z,\delta)-\underline U(\hat z,\delta)\to 0$ as $\delta\to 0$. 

Using  the renewal theory, we obtain, if the law of $\xi$ is  non-arithmetic,  that
\begin{equation}\label{ren_theory_non_arith}
\lim_{x\to\infty}\,e^{\beta x}\bar{H}(x)=\frac{1}{\tilde \E \xi}\,\int^{\infty}_{0}\,z(v)\d v,
\end{equation} see, e.g.,  Ch. XI, 9, \cite{Feller1971}.
If the law of $\xi$ is arithmetic with the step $d>0$, then, according to Proposition \ref{Iks} for any $x>0$
\begin{equation}\label{ren_theory_arith}
\lim_{n\to\infty}\,e^{\beta (x+nd)}\bar{H}(x+nd)=\frac{d}{\tilde \E \xi}\sum_{j\in {\mathbb Z}}z(x+jd)\,
I_{\{x+jd\ge 0\}}.
\end{equation}
The equalities \eqref{ren_theory_non_arith}
and
\eqref{ren_theory_arith}
implies the statement. 
\fdem

\smallskip 
 The proof of the result below, formulated to cover our needs, follows the same line as in Lemma 2.6  of the Buraczewski--Damek paper \cite{BD} with minor changes to include also the arithmetic case.  

\begin{theo} 
\label{Le.3.2+}
Suppose that (\ref{3.3}) hold. 
If the  support of distribution of $Y_{\infty} $ is unbounded from above then
$$
\liminf_{u\to\infty}\,u^{\beta}\,\P(Y_{\infty}>u)>0\,.
$$
\end{theo} 
\noindent{\sl Proof.} 
Let 
$$
\bar Y_n:=-\sum_{j=1}^n\,Q_j^-\,Z_{j-1}, \qquad Y_{n,\infty}:=\sum^\infty_{j=n+1}\,Q_j\,\prod^{j-1}_{l=n+1}\,M_l
$$
and let $Z^*_n:=\sup_{j\le n}Z_j$.  
Theorems \ref{Le.3.2}, \ref{Le.3.2++} imply that $ \P(\bar Y_{\infty} <-u)\le C_1u^{-\beta}$ with $C_1>0$. for sufficiently large $u$.
On the other hand,  by Proposition
\ref{Pr.prod_M_j}
 $\P(Z^*_{\infty} >u)\ge  C_2u^{-\beta}$ with $C_2>0$ and $u\to\infty$.
\smallskip 
%
%

Put $U_n:=\{Z_n>u,\, \bar Y_n> -Cu\}$   where $C^\beta:=4C_1/C_2$. The process  $\bar Y$  decreases.  Therefore, we have the inclusion $ \{Z_n>u\}\subseteq \{\bar Y_\infty \le -Cu \}\cup U_n$. It follows that  
for sufficiently large $u>0$
\bean
(3/4)C_2u^{-\beta}&\le& \P(Z^*_\infty>u)=\P(\cup_n \{Z_n>u\})
\le \P(\bar Y_\infty \le -Cu)+\P(\cup_n U_n)\\
&\le& 2C_1C^{-\beta} u^{-\beta}+\P(\cup_n U_n) 
\eean 
implying that $
\P(\cup_n U_n)\ge (1/4)C_2u^{-\beta}$. 
\smallskip

Since $\bar Y_n+Z_nY_{n,\infty}\le Y_n+Z_nY_{n,\infty}=Y_{\infty}$, we have that 
$$
\{Y_{n,\infty}>C+1\}\cap U_n\subseteq \{\bar Y_n+Z_nY_{n,\infty}>u\}\cap U_n\subseteq \{Y_{\infty}>u\}\cap U_n,
$$
Note  that $\P(Y_{\infty}>C+1)=\P(Y_{n,\infty}>C+1)$, because $\cL(Y_{n,\infty}) =\cL(Y_{\infty})$.  
Using  the independence  of $Y_{n,\infty}$ and the sets $W_n:=U_n\cap\big(\cup_{k=1}^{n-1}U_k\big)^c$
forming a disjoint partition of $\cup_nU_n$, we get that  
 \bean
\P(Y_{\infty}>C+1)\P(\cup_n W_n)
&=&\sum_n\P\big(\{Y_{n,\infty}>C+1\}\cap W_n\big)\\
&\le&  \sum_n\P\big(\{Y_{\infty}>u\}\cap W_n)
\le  \P(Y_{\infty}>u). 
\eean
Thus,  $\P(Y_{\infty}>u)\ge (1/4)bC_2u^{-\beta}$ where $b:=\P(Y_{\infty}>C+1)>0$  by the assumption that the support of $\cL(Y_{\infty})$ is unbounded from above. The obtained asymptotic bound implies that   $C_+>0$. \fdem

\medskip
Summarizing  the above  results we get for function $\bar G(u)=\P(Y_{\infty}>u)$  the following asymptotic properties when $u\to \infty$: 
\begin{theo}
\label{end}
Suppose that (\ref{3.3}) holds. Then $\limsup  u^\beta \bar G(u)<\infty$.  If  $Y_{\infty}$ is unbounded from above, then  
$\liminf  u^\beta \bar G(u)>0$ and  in the case where  $\cL(\ln M)$ is non-arithmetic 
$\bar G(u)\sim C_+u^{-\beta}$ where $C_+>0$. 
 \end{theo}

\medskip

{\bf Acknowledgements.}  
The research  is funded by the grant of the Government of Russian Federation $n^\circ$14.A12.31.0007. The second author is
partially supported by the grant of RSF number 14-49-00079, National Research University ``MPEI", 
14 Krasnokazarmennaya, 111250 Moscow, Russia.


\begin{thebibliography}{100}

\bibitem{ABL} 
Albrecher H., Badescu A., Landriault D. On the dual risk model with taxation,
Insurance: Mathematics and Economics, {\bf 42} (2008), 1086--1094.

\bibitem{Asm}
 Asmussen S., Albrecher H.  {\em Ruin Probabilities.} World Scientific,
Singapore, 2010.

\bibitem{AGS} 
Avanzi B., Gerber H.U., Shiu E.S.W. Optimal dividends in the dual model.
{\sl Insurance: Mathematics and Economics}, {\bf 41} (2007), 111--123.

\bibitem{BKM} 
Bankovsky D., Kl\"uppelberg C., Maller R. On the ruin probability of the generalised Ornstein--Uhlenbeck process in the Cram\'er case. {\sl Journal of Applied Probability}, {\bf 48A} (2011), 15--28.  

\bibitem{BD} 
Buraczewski D., Damek E. A simple proof of heavy tail estimates for affine type Lipschitz recursions. {\sl Stoch. Proc. Appl.}, {\bf 127} (2017), 657--668. 

\bibitem{BuraczewskiDamekMikosch_2016} 
Buraczewski D., Damek E.,  Mikosch Th.  
{\em Stochastic Models with Power-Law Tails. The Equation $X = AX + B$}.
Springer Series in Operations Research and Financial Engineering,
Springer,  2016.

 \bibitem{Bayraktar-Egami} 
Bayraktar E., Egami M. Optimizing venture capital investments in a jump diffusion model. 
{\sl Mathematical Methods of Operations Research}, {\bf 67} (2008), 1, 21--42.


\bibitem{BJ} 
Bichteler K., Jacod J.  Calcul de Malliavin pour les diffusions avec sauts: existence d'une densit\'e dans le cas unidimensionnel. {\em S\'eminaire de probabilit\'e, XVII}, Lecture Notes in Math., {\bf  986}, Springer, Berlin, 1983, 132--157.


\bibitem{CT}
Cont R., Tankov P. {\em Financial Modelling with Jump Processes.} 
Chapman \& Hall, 2004. 


\bibitem{Feller1971}
Feller W. {\em An Introduction to Probability Theory and Its Applications.}
{\bf 2}, 2nd edition, Wiley, New York, 1971.


\bibitem{FrKP} 
Frolova A., Kabanov Yu., Pergamenshchikov S. In the  insurance business risky investments are dangerous, {\em Finance and Stochastics}, {\bf 6}  (2002),  227--235.

\bibitem{Go91} Goldie C.M.  
Implicit renewal theory and tails of solutions of random equations.
{\em The Annals of Applied Probability}, {\bf 1} (1991), 1, 126--166.


\bibitem{Gr} Grandell I. 
{\em Aspects of Risk theory.} Springer,
Berlin, 1990.

\bibitem{Grincevicius1975_a}
Grincevicius  A.K.  One limit theorem for a random walk on the line.
 {\em  Institute of Physics and Mathematics, Academy of Sciences of the Lithuanian SSR},  {\bf 15} (1975),  4,  79--91. 

\bibitem{GL}
Guivarc'h Y., Le Page E. On the homogeneity at infinity of the
stationary probability for affine random walk. In:  Bhattacharya S., Das T., Ghosh A.,  Shah R. (eds). {\em Recent trends in ergodic theory and dynamical systems}, Contemporary Mathematics, 119--130,  AMS, 2015.

\bibitem{IksanovPolotskiy2016}
Iksanov A.,  Polotskiy S. 
Tail behavior of supreme of perturbed random walks.
{\em Theory of Stochastic Processes}, 
{\bf  21 (37)}  (2016), 1, 12--16.

 \bibitem{JS} Jacod J., Shiryaev A.N. {\em Limit theorems for stochastic processes.} 2nd edition, Springer, Berlin, 2002.
 
 \bibitem{KP} Kabanov Yu., Pergamenshchikov S. In the  insurance business risky investments are dangerous: the case of negative risk sums. {\em Finance and Stochastics},   {\bf 20} (2016), 2,  355 -- 379. 

\bibitem{Kalash-Nor} 
Kalashnikov V., Norberg R. Power tailed ruin
probabilities in the presence of risky investments. {\em Stoch. Proc. Appl.}, {\bf 98} (2002), 211--228.

\bibitem{Kesten}
Kesten H. Random difference equations and renewal theory for products
of random matrices. {\em Acta Math.} {\bf 131} (1973) 207--248. 

\bibitem{KKM}
Kl\"uppelberg C,  Kyprianou A.E., Maller R.A. 
Ruin probabilities and overshoots for general L\'evy insurance risk processes.  {\sl Ann. App. Probab.}, {\bf  14} (2004), 4,  1766--1801.



\bibitem{LL}
Lamberton D.,  Lapeyre B. {\sl Introduction to Stochastic Calculus Applied to Finance.}  Chapman \& Hall, London, 1996. 

\bibitem{MR}
Marinelli C., R\"ockner M. 
On maximal inequalities for purely discontinuous martingales in infinite dimensions. 
{\em S\'eminaire de Probabilit\'es}, Lect. Notes Math., {\bf XLVI} (2014), 293--315. 

\bibitem{Novikov} 
Novikov A.A. On discontinuous martingales.  {\em Theory Probab. Appl.}, {\bf 20}, (1975),
1, 11--26. 

\bibitem{Nyrh-99} 
Nyrhinen H. On the ruin probabilities in a
general economic environment. {\em Stoch. Proc. Appl.}, {\bf 83} (1999), 319--330.

\bibitem{Nyrh-01} 
Nyrhinen H. Finite and infinite time ruin probabilities in a stochastic economic environment. {\em Stoch. Proc. Appl.}, {\bf 92} (2001), 265--285.

\bibitem{Paul-93} 
Paulsen J. Risk theory in stochastic economic environment.  {\em Stoch. Proc. Appl.}, {\bf  46} (1993), 327--361.


\bibitem{Paul-98}
 Paulsen J. Sharp conditions for certain ruin in a risk process with stochastic return on investments. {\em Stoch. Proc. Appl.}, {\bf 75} (1998), 135--148.

\bibitem{Paul-02}
Paulsen J.  On Cram\'er-like asymptotics for risk processes with stochastic return on investments. {\em Ann. Appl. Probab.},  {\bf 12} (2002), 4, 1247--1260.

\bibitem{Paulsen1993}
Paulsen J., Gjessing H. K. Ruin theory with stochastic return on
investments. {\em Adv. Appl. Probab.},  {\bf 29} (1997), 4, 965--985.

\bibitem{PeZe-06} 
Pergamenshchikov S., Zeitouny O.  
Ruin probability in the presence of risky investments. {\em Stoch. Process. Appl.}, {\bf  116} (2006), 267--278. Erratum to: ``Ruin probability in the presence of risky
investments". {\em  Stoch.
Proc. Appl.},  {\bf 119} (2009), 1, 305--306.

\bibitem{Sato1999}
Sato K. {\em L\'evy processes and Infinitely Divisible Distributions}.
Cambridge University Press, Cambridge,  1999. 

\bibitem{Sax} 
Sax\'en T.  On the probability of ruin in the collective risk theory for
insurance enterprises with only negative risk sums. {\em Scand. Actuarial J.} 1948, 1-2,  199--228. 

\end{thebibliography}
\end{document}